\documentclass[onefignum,onetabnum]{siamart220329}

\usepackage{lipsum}
\usepackage{amsfonts}
\usepackage{graphicx}
\usepackage{epstopdf}
\usepackage{booktabs}
\usepackage{multirow}
\usepackage{makecell}
\usepackage{algorithmic}
\usepackage{comment}

\ifpdf
  \DeclareGraphicsExtensions{.eps,.pdf,.png,.jpg}
\else
  \DeclareGraphicsExtensions{.eps}
\fi

\definecolor{red2}{RGB}{204,0,0}
\definecolor{blue2}{RGB}{0,103,165}
\hypersetup{colorlinks,linkcolor=[RGB]{0,103,165},citecolor=[RGB]{180,0,0}}
\usepackage[nocompress]{cite}

\newsiamremark{remark}{Remark}
\newsiamremark{hypothesis}{Hypothesis}
\crefname{hypothesis}{Hypothesis}{Hypotheses}
\newsiamthm{claim}{Claim}

\def\be{\begin{equation}}
\def\ee{\end{equation}}

\def\x{\mathbf{x}}

\def\Xt{\widetilde{\X}}

\def\y{\mathbf{y}}

\def\f{\mathbf{f}}
\def\g{\mathbf{g}}
\def\z{\mathbf{z}}

\def\X{\mathbf{X}}
\def\Z{\mathbf{Z}}
\def\Y{\mathbf{Y}}

\def\Rs{\mathbb{R}}

\def\G{\mathbf{G}}
\def\E{\mathbf{E}}

\def\N{\mathbf{N}}

\def\T{\mathbf{T}}
\def\Pi{\mathbf{\Phi}}

\def\m{\boldsymbol{\mu}}
\def\s{\boldsymbol{\sigma}}
\def\xxi{\boldsymbol{\xi}}
\def\D{\mathbf{D}}
\def\S{\mathbf{S}}

\newcommand{\norm}[1]{\left\lVert#1\right\rVert}

\newcommand{\met}[1]{\textcolor{black}{#1}}

\usepackage{amsopn}

\headers{Modeling Stochastic Chemical Reaction Network}{Yuan Chen Weize Mao and Dongbin Xiu}

\title{Data-driven Effective Modeling of Stochastic Chemical Reaction Networks}

\author{Yuan Chen\thanks{E-mail addresses: \texttt{chen.11050@osu.edu}, \texttt{weizemao339@gmail.com}, \texttt{xiu.16@osu.edu}. Department of Mathematics, The Ohio State University, Columbus, OH 43210, USA. Funding: This work was partially supported by AFOSR FA9550-24-1-0237.} \and Weize Mao\footnotemark[1] \and Dongbin Xiu\footnotemark[1]}

\begin{document}

\maketitle

\begin{abstract}
The Stochastic Simulation Algorithm (SSA), widely considered an exact algorithm for stochastic chemical reaction networks, suffers from high computational cost. In this work, we propose a data-driven effective model that operates on a user-defined coarse time step independent of the underlying microscopic reaction-event scale. This is accomplished by directly approximating the finite-time transition kernel of the continuous-time Markov chain induced by SSA, using a generative machine learning model trained on short bursts of SSA simulation data. The trained model constructs a stochastic propagator that recursively generates statistically consistent trajectories at the constant coarse time step, with significantly reduced computational cost. In this paper, we employ conditional normalizing flow as the stochastic propagator. A comprehensive set of numerical examples is presented to demonstrate the accuracy and efficiency of the proposed method. 
\end{abstract}



\section{Introduction}

Stochastic chemical kinetics of well-mixed reacting systems are commonly modeled by
the Chemical Master Equation (CME), which describes the evolution of the probability
distribution of the underlying reaction network. The corresponding stochastic process
is a continuous-time Markov jump process whose exact sample trajectories can be generated
by the Stochastic Simulation Algorithm (SSA), also known as the Gillespie algorithm.
Since its introduction in \cite{gillespie1977exact}, SSA has become one of the most
widely used computational tools for stochastic chemical reaction systems.

Despite its exactness, SSA suffers from high computational cost. The algorithm advances the system by simulating one reaction event at a time, and the corresponding firing times can become exceedingly small for systems with large populations or stiff reaction rates. Consequently, long-time simulations often require an extremely large number of reaction events. Numerous acceleration and approximation techniques have been developed to reduce this cost, including the next-reaction method \cite{gibson2000efficient}, modified direct method \cite{cao2004efficient}, sorting direct method \cite{mccollum2006sorting}, moment closure approximations \cite{grima2012study}, and tau-leaping methods. Among these approaches, tau-leaping advances the system over a finite time interval by approximating the collective statistics of many reaction events \cite{gillespie2001approximate}. This idea has led to extensive developments, including stability and error analysis \cite{rathinam2003stiffness,li2007analysis,anderson2011error}, nonnegativity-preserving variants \cite{chatterjee2005binomial,cao2005avoiding}, and adaptive or implicit schemes for stiff systems \cite{cao2007adaptive}. Related coarse-graining and reduction strategies for stochastic reaction systems have also been studied from both computational and mathematical perspectives \cite{katsoulakis2003coarse,katsoulakis2019data}. We refer to \cite{higham2008modeling,anderson2015stochastic,schnoerr2017approximation} for broader reviews and background. These methods are effective in many regimes, but they typically rely on assumptions regarding slowly varying propensity functions and can suffer from loss of accuracy or stability when these assumptions fail.

In this work, we pursue a different direction. Instead of approximating the individual
reaction events or modifying the internal mechanics of SSA, we directly approximate
the finite-time transition law of the Markov jump process induced by SSA. More precisely,
let $X(t)$ denote the stochastic process generated by the reaction network. For any
finite time lag $\Delta>0$, the process admits a transition kernel
\begin{equation}
P_{\Delta}(x,x')
=
\mathbb P(X(t+\Delta)=x' \mid X(t)=x),
\end{equation}
which defines the finite-time Markov semigroup associated with the CME.
The objective of this work is to construct a data-driven approximation of
$P_{\Delta}$ directly from short bursts of SSA simulation data.

This viewpoint leads naturally to a coarse stochastic propagator operating on a user-defined time step $\Delta$, which can be substantially larger than the average SSA firing time.
Unlike tau-leaping methods, which accelerate SSA by approximating the numbers of reaction firings over a finite interval under assumptions on the variation of propensity functions, the proposed approach does not approximate the reaction-event dynamics. Instead, it directly targets the finite-time transition kernel of the SSA-induced continuous-time Markov chain and learns this kernel from SSA data. Consequently, the learned model should be viewed as a data-driven numerical approximation of the exact Markov semigroup associated with the underlying continuous-time Markov chain, rather than as an approximation of the reaction-event dynamics.

The proposed method follows the framework of stochastic Flow Map Learning (sFML)
\cite{chen2024learning}, which extends deterministic flow map learning
\cite{qin2019data} to stochastic dynamical systems. In deterministic systems, flow
map learning seeks to approximate the map between consecutive states over a finite
time interval. In stochastic systems, however, the finite-time evolution is no longer
deterministic and is instead characterized by a conditional probability distribution.
Consequently, the stochastic flow map is interpreted as a stochastic realization operator
whose generated samples follow the finite-time transition kernel.

To realize the stochastic flow map, we employ a conditional normalizing flow model.
Normalizing flows (see, e.g., \cite{papamakarios2021normalizing}) provide expressive generative models capable of representing highly
non-Gaussian conditional distributions while maintaining tractable likelihood evaluation
and efficient sampling. This property is particularly important for stochastic chemical
reaction systems, whose finite-time transition distributions are frequently strongly
non-Gaussian and highly state-dependent.

The resulting framework constructs a stochastic propagator directly from
short bursts of SSA trajectory data. Once trained, the model can efficiently produce
long-time trajectories at a prescribed coarse time step while preserving important
statistical properties of the underlying stochastic dynamics.
The main contributions of this work are summarized as follows:
\begin{itemize}
\item We reinterpret stochastic flow map learning for chemical reaction systems as
a data-driven approximation of the finite-time transition kernel associated with the
SSA Markov process.

\item We construct a coarse stochastic propagator operating on a user-defined time
step independent of the microscopic reaction-event scale.

\item We employ conditional normalizing flows to approximate non-Gaussian
transition distributions arising from stochastic chemical kinetics.

\item We demonstrate through a comprehensive set of numerical examples that the
learned propagator produces statistically consistent trajectories while significantly
reducing the number of simulation steps compared with SSA.
\end{itemize}

\section{Preliminary} \label{sec:setup}
We first present the problem setup and our main objective.

\subsection{Problem Setup}

Consider a stochastic chemical reaction network with $N$ species
$S_i$, $i=1,\dots,N$, and $M$ reaction channels
\begin{equation}
R_{\mu}:
\qquad
\sum_{i=1}^{N} r_{i\mu} S_i
\overset{c_{\mu}}{\longrightarrow}
\sum_{i=1}^{N} p_{i\mu} S_i,
\qquad
\mu=1,\dots,M,
\end{equation}
where $c_{\mu}$ denotes the reaction rate constant, and
$r_{i\mu},p_{i\mu}$ are nonnegative integers representing the stoichiometric
coefficients.

Let $N$-dimensional vector
\[
\X(t)=[X_1(t),\dots,X_N(t)]^T
\]
denote the vector of molecule numbers at time $t$, in the state space $\mathbb{Z}^N_{\geq 0}$, the set of $N$-dimensional vectors of non-negative integers. Under the well-mixed assumption, the evolution of $\X(t)$ is modeled as a continuous-time Markov chain (CTMC). Each
reaction channel $R_{\mu}$ is characterized by:
\begin{itemize}
\item the stoichiometric vector
$
\mathbf{s}_{\mu}
=
[p_{1\mu}-r_{1\mu},\dots,p_{N\mu}-r_{N\mu}]^T,
$
which specifies the state change caused by the reaction;

\item the propensity function
$
a_{\mu}(\x)=c_{\mu}h_{\mu}(\x),
$
where $h_{\mu}(\x)$ represents the number of available reactant combinations at state $\X(t)=\x$.
\end{itemize}
The stochastic process admits an infinitesimal generator
\begin{equation}
(\mathcal L f)(\x)
=
\sum_{\mu=1}^{M}
a_{\mu}(\x)\left(f(\x+\mathbf{s}_{\mu})-f(\x)\right),
\end{equation}
which defines the corresponding chemical master equation.

The SSA algorithm generates exact sample trajectories of this CTMC.
Starting from state $\X(t)=\x$, SSA samples:
\begin{enumerate}
\item[(1)] a firing time
\[
\tau
\sim
\text{Exp}\left(
\sum_{\mu=1}^{M} a_{\mu}(\x)
\right),
\]
\item[(2)] a reaction index $\lambda$ with probabilities
\[
\mathbb P(\lambda=\mu)
=
\frac{a_{\mu}(\x)}
{\sum_{\nu=1}^{M} a_{\nu}(\x)},
\]
\end{enumerate}
and updates the state according to
\[
\X(t+\tau)=\X(t)+\mathbf{s}_{\lambda}.
\]

Although SSA evolves through random microscopic firing events,
the resulting process is globally characterized by a finite-time transition kernel.
For any $\Delta>0$, define transition kernel $P_\Delta: \mathbb{Z}^N_{\geq 0} \times \mathbb{Z}^N_{\geq 0} \to [0,1]$
\begin{equation} \label{P_Delta}
P_{\Delta}(\x,\x')
=
\mathbb P(\X(t+\Delta)=\x'|\X(t)=\x).
\end{equation}
Equivalently, for any observable $\phi$,
\begin{equation}
(\mathcal P_{\Delta}\phi)(\x)
=
\mathbb E[\phi(\X(t+\Delta))|\X(t)=\x]
\end{equation}
defines the Markov semigroup associated with the process. 

For each fixed $\x$, $P_\Delta(\x,\cdot)$ defines a (conditional) probability measure on $\mathbb{Z}^N_{\geq 0}$ where $P_\Delta(\x,A)=\mathbb{P}(\X(t+\Delta)\in A|\X(t)=\x)$ for any event set $A$. Modeling this conditional distribution is the central goal of this paper.

\subsection{Objective}

The objective of this work is to construct a data-driven approximation of the
finite-time transition kernel $P_{\Delta}$ \eqref{P_Delta} for a prescribed coarse time step
$\Delta$. More importantly, $\Delta$ is chosen according to the time scale of the macroscopic behavior of the system and is independent of the microscopic SSA firing 
times $\tau$. Therefore, $\Delta$ may be substantially larger than the average reaction-event scale, i.e., $\Delta \gg \mathbb{E}[\tau]$.

\subsection{Related Work}
This work belongs to the broader area of learning dynamical systems from data, which has attracted growing attention due to the rapid development of machine learning. Most existing methods focus on deterministic systems, see, for example, \cite{raissi2019physics,raissi2018multistep,brunton2016discovering,li2020fourier,owhadi2021computational,qin2019data,Churchill_2023,lu2021learning,chen2026targeted}.

For stochastic systems, most existing approaches are `model-based', i.e., observational data is assumed to be sampled from a known model with unknown components, such as It\^{o} stochastic differential equations which approximately admit Gaussian transition distributions with unknown drift and diffusion functions. These unknown functions are approximated by analyzing and matching statistical properties of data, such as moments, likelihoods, and densities. Such methods have been developed in conjunction with Gaussian processes \cite{yildiz2018learning,pmlr-v1-archambeau07a,darcy2022one,opper2019variational} and deep neural networks (DNNs) \cite{chen2023data,yang2022generative,chen2021solving,zhang2022multiauto,dietrich2023learning,xia2024efficient,lu2024learning}.

For processes generated by SSA, and CTMCs more broadly, such methods are less capable since the transition distributions are \emph{non-Gaussian} and the underlying mathematical structure does not conform to a standard parametric form. Some recent explorations using machine learning methods can be found in \cite{tang2023neural,sukys2022approximating,cairoli2023generative}. In this work, we adopt the sFML framework \cite{chen2024learning}, which incorporates generative models to approximate the stochastic map between consecutive states in distribution. Many generative models are compatible with this framework, including GANs (generative adversarial networks) \cite{chen2024learning}, autoencoders \cite{xu2023learning}, normalizing flows \cite{chen2026modeling,chenxiu2024}, and diffusion models \cite{liu2024training,song2021score}. Recently, it was also extended to stochastic systems subjected to external excitation and multiscale systems \cite{chen2026modeling,chenxiu2024}.

\section{Main Method: Effective Stochastic Propagator}
\label{method}

\subsection{Finite-Time Stochastic Flow Map}
Let $t_n=n\Delta$, $n=0,1,\cdots$ be a sequence of discrete time points with constant time stepsize $\Delta$ and we denote $\X(t_n)$ by $\X_n$ for convenience. Then the transition from $\X_n$ to $\X_{n+1}$ is governed by the transition kernel
\[
\X_{n+1}\sim P_{\Delta}(\X_n,\cdot).
\]

We therefore define the stochastic flow map as a stochastic realization operator
\begin{equation}
 \X_{n+1}=\G_{\Delta}(\X_n,\xxi_n),
\end{equation}
where $\xxi_n$ is an auxiliary random variable independent of $\X_n$ and follows a prescribed distribution $\mu_{\xxi}$. The stochastic flow map is required to match the conditional distribution
\begin{equation} \label{app}
\G_{\Delta}(\x,\xxi)  \sim P_{\Delta}(\x,\cdot),\quad \text{as }\xxi \sim \mu_{\xxi}.
\end{equation} 
In other words, $\G_{\Delta}(\x,\xxi)$ provides a sampling realization of the transition kernel $P_{\Delta}(\x,\cdot)$.

In the following sections, we construct a numerical approximation $\G_{\Theta}$ of the stochastic flow map, parameterized by $\Theta$, such that the conditional distribution of $\G_{\Theta}(\x,\xxi)$ approximates the finite-time transition kernel $P_\Delta(\x,\cdot)$. Consequently, the learned model acts as a sampling realization of the finite-time
transition kernel of the SSA Markov process. Once trained, repeated application of the stochastic flow map produces a discrete-time stochastic process
\begin{equation} \label{sFM1}
\widetilde \X_{n+1} = \G_{\Theta}(\widetilde \X_n,\xxi_n),
\end{equation}
whose finite-time statistics approximate those of the original SSA process sampled
at time interval $\Delta$.

\subsection{Short-Burst Data Generation}
\label{data}

The training data are generated directly from SSA simulations.
Let $\Delta>0$ denote the prescribed coarse time step and let
\[
\left\{\X_0^{(i)}\right\}_{i=1}^{K}
\]
be randomly sampled initial conditions from a prescribed distribution with size of dataset $K>0$.

For each initial condition $\X_0^{(i)}$, we run SSA up to time $\Delta$ and obtain
\[
\X_1^{(i)}=\X^{(i)}(\Delta).
\]
This produces a collection of independent sample pairs
\begin{equation} \label{dataset}
\left(
\X_0^{(i)}, \X_1^{(i)}
\right),
\qquad
i=1,\dots,K.
\end{equation}

Each sample pair represents one realization of the finite-time transition kernel
\[
\X_1^{(i)}
\sim
P_{\Delta}\left(\X_0^{(i)},\cdot\right).
\]

An important feature of this construction is that only short independent trajectory
bursts are required. The training process therefore avoids the need for extremely
long SSA simulations and is naturally parallelizable.

\subsection{Conditional Normalizing Flow Approximation}

To approximate the stochastic flow map, we employ a conditional normalizing flow model.
Let
\[
\z\sim \mathcal N(0,\mathbf{I}_N)
\]
be an $N$-dimensional standard Gaussian random variable.

We seek a parameterized map
\begin{equation}
\X_1 = \G_{\Theta}(\X_0,\z),
\end{equation}
such that the conditional distribution of $\G_{\Theta}(\X_0,\z)$ approximates
the finite-time transition kernel
\[
P_{\Delta}(\X_0,\cdot).
\]

In this work, the map $\G_{\Theta}$ is realized through a conditional normalizing
flow. More specifically, we introduce an invertible transformation
\[
\T_{\theta}: \mathbb R^N \to \mathbb R^N,
\]
whose parameters depend on the conditioning variable $\X_0$ through a
neural network with hyper-parameter $\Theta$:
\begin{equation}
\theta = \mathbf N_{\Theta}(\X_0).
\end{equation}
The stochastic flow map is then defined as
\begin{equation} \label{x1}
\X_1 = \T_{\mathbf N_{\Theta}(X_0)}(\z).
\end{equation}

The invertibility of the map $\mathbf{T}$ enables tractable likelihood evaluation through the change-of-variable formula. Let
\[
\S_{\theta}=\T_{\theta}^{-1}.
\]
Then the conditional density is
\begin{equation}
p(\X_1|\X_0;\Theta)  = p_{\z}\left( \S_{\mathbf N_{\Theta}(\X_0)}(\X_1) \right)
\left|  \det \D\T_{\mathbf N_{\Theta}(\X_0)} \right|^{-1},
\end{equation}
where $p_{\z}$ denotes the standard Gaussian density function and $\D\T(\y) = \partial\T/\partial \y$ denotes the Jacobian of the mapping $\T$.

\begin{figure}[htbp]
  \centering
  \includegraphics[width=.8\textwidth]{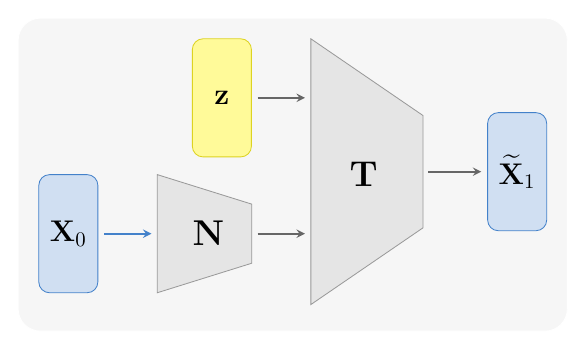}
  \caption{The DNN model structure for the proposed normalizing flow
    sFML method \eqref{x1}.}
  \label{fig:NFnet}
\end{figure}

The model parameters $\Theta$ are obtained by maximizing the conditional likelihood
over the training dataset, equivalently minimizing the empirical negative log-likelihood loss
\begin{equation} \label{loss}
\mathcal L(\Theta) = - \frac{1}{K}\sum_{i=1}^{K} \log p\left(\X_1^{(i)}|\X_0^{(i)};\Theta\right).
\end{equation}
An illustration of the proposed sFML model structure can be found in Figure
\ref{fig:NFnet}. This is in direct correspondence of \eqref{x1}.

We emphasize that the conditional normalizing flow serves only as one realization
strategy for the stochastic flow map. The proposed framework itself is independent
of the particular generative architecture and may incorporate alternative conditional
generative models.

\subsection{Recursive System Prediction}

Once the model is trained, long-time trajectories are generated recursively by
\be \label{predict}
\left\{
\begin{split}
&\Xt_0 = \X_0, \\
   & \widetilde \X_{n+1} = \mathbf{R}\left( \G_{\Theta}(\widetilde \X_n,\z_n)\right), \qquad n=0,1,2,...,
\end{split}
\right.
\ee
where $\z_n$ are standard Gaussian random samples independent of $\widetilde \X_n$ and
$\mathbf{R}$ denotes an enforcing operator to satisfy the  physical constraints such
as nonnegativity and integer-valued molecule counts -- see the Section \ref{compdetail} below.

The resulting process is designed to provide a weak approximation to the true stochastic dynamics at the discrete time points, i.e., 
$$
(\Xt_0,\Xt_1,\cdots,\Xt_n) \stackrel{d}{\approx} (\X_0,\X_1,\cdots,\X_n),
$$ 
but operates only on a constant coarse time step $\Delta$ independent of the microscopic SSA firing times. Consequently, the learned stochastic propagator can significantly reduce the number of simulation steps required for long-time stochastic simulations.


\subsection{Computational Details for Enforcing Operator}
\label{compdetail}

\subsubsection*{Conservation of Molecule Numbers}
In some reaction networks, the total molecule count is a conserved quantity, i.e.,
\begin{equation}
    \X_n \cdot \mathbf{1}_N = \X_0 \cdot \mathbf{1}_N, \text{~~~for all } n = 0, 1, \cdots,
\end{equation}
where $\mathbf{1}_N$ denotes the $N$-dimensional vector of all ones. In this case, the transition kernel $P_\Delta(\x, \cdot)$ is supported on the $(N-1)$-dimensional hyperplane $\{\x': \x' \cdot \mathbf{1}_N = \x \cdot \mathbf{1}_N\}$, a constraint that will not be preserved automatically by generic generative models. We enforce this constraint strongly into the model construction.

Let $\mathbf{A} \in \mathbb{R}^{(N-1)\times N}$ be a slicing matrix consisting of the first $N-1$ rows of the $N\times N$ identity matrix. We replace the training dataset with the projected pairs
\begin{equation}\label{equ:data2o_m}
    \left(\Y_0^{(i)},\, \Y_{1}^{(i)}\right) := \left(\X_0^{(i)},\, \mathbf{A}\X_{1}^{(i)}\right), \quad i = 1, \ldots, K,
\end{equation}
and train an $(N-1)$-dimensional conditional normalizing flow model based on the new dataset
\begin{equation} \label{reducednfmodel}
    \Y_1=\widehat{\T}_{\N_\Theta(\Y_0)}\left(\widehat{\z}\right)
\end{equation} 
following the procedure stated in previous sections. Here $\widehat{\z}$ follows an $(N-1)$-dimensional unit Gaussian distribution.

After training, the full $N$-dimensional model prediction is then recovered using the dimension-extension operator $\E : \mathbb{R}^N \times \mathbb{R}^{N-1} \to \mathbb{R}^N$, defined by
\begin{equation}
    \E(\x, \y) := \left[\y;\; \x\cdot\mathbf{1}_N - \y\cdot\mathbf{1}_{N-1}\right]^T,
\end{equation}
which reconstructs the $N$-th component so that the total count is exactly preserved. The resulting model takes the form
\begin{equation}
    \G_\Theta(\x, \widehat{\z}) = \E\!\left(\x,\, \widehat{\T}_{\N_\Theta(\x)}\left(\widehat{\z}\right)\right),
\end{equation}
which satisfies $\G_\Theta(\x, \widehat{\z}) \cdot \mathbf{1}_N = \x \cdot \mathbf{1}_N$ for any $\Theta$ by construction. Other invariants of this type, including those expressible as signal temporal logic (STL) requirements~\cite{cairoli2023generative}, are left for future work.

\subsubsection*{Post-processing}
A notable feature of chemical reaction systems is that the state $\X(t)$ takes values of non-negative integers, as each component represents a molecule count. Since the output of the normalizing flow lies in $\mathbb{R}^N$, we apply a deterministic post-processing operator $\mathbf{R} : \mathbb{R}^N \to \mathbb{Z}_{\geq 0}^N$ after each recursive step, defined by
\begin{equation}
    \mathbf{R}(\x) = \mathrm{Round}\!\left(\max(\x, \mathbf{0})\right),
\end{equation}
where $\max(\cdot, \mathbf{0})$ enforces nonnegativity componentwise and $\mathrm{Round}(\cdot)$ maps each component to the nearest integer. The operator $\mathbf{R}$ ensures that all generated states remain physically admissible. This projection onto the discrete state space may be interpreted as an abstraction of the learned real-valued model, we refer to \cite{bortolussi2018deep} for interested readers.

\section{Numerical Examples}
\label{sec:examples}
In this section, we present several numerical experiments to demonstrate the performance of the proposed sFML model across the following reaction networks:
\begin{itemize}
\item Transfer process;
\item Lotka-Volterra model;
\item The Brusselator;
\item The Autocatalysis;
\item The Oregonator.
\end{itemize}

In all examples, the reaction networks are known, and the training data are generated via SSA following the procedure described in Section \ref{data}. Initial conditions are sampled uniformly from a bounded domain specified for each example, along with the chosen coarse time step $\Delta$. Once trained, the model is used to generate ensemble predictions with new initial conditions, which are compared against the SSA ground truth \met{and the classical acceleration method tau-leaping} using the following quantities:
\begin{itemize}
    \item \met{Number of simulation steps required to reach the same termination time, used as a measure of computational efficiency;}
    \item Simulated trajectories and their mean and standard deviation, for visual comparisons as well as basic statistics;
    \item Fourier modes for examples exhibiting periodic behavior;
    \item Distributional comparison of the learned stochastic flow map
     $\G_\Theta(\x,\cdot)$ and the SSA transition kernel $P_{\Delta}(\x,\cdot)$ at fixed conditioning states $\x$.
\end{itemize}

We use the Masked Autoregressive Flow (MAF, c.f. \cite{papamakarios2017masked}) model as a technical realization of the normalizing flow. Upon use of the model, there are $4$ coupling flow mappings constructed with independent training parameters. The training data size, batch size, and number of training epochs will be specified in each example. We use fully connected neural networks with 3 layers, each of which has 20 nodes, with a $\tanh$ activation function. We employ the cyclic learning rate with a base rate $3\times10^{-4}$ and a maximum rate $5\times10^{-4}$, $\gamma=0.99999$, and step size $10,000$. Simultaneously, the whole cycle is set to decay with a scale $0.5$ for every $40,000$ training epoch. To reduce the overfitting phenomenon, we also add a small weight decay with a scale $0.01$ on the updated gradients.

\met{For SSA, one simulation step corresponds to one reaction firing event. For sFML, one simulation step corresponds to one recursive evaluation of the learned stochastic flow map $\G_\Theta$ over the prescribed coarse time step $\Delta$, and therefore the number of sFML steps is exactly $T/\Delta$. We also include tau-leaping as a classical accelerated SSA method. When the propensity functions remain approximately unchanged over a time interval $\tau$, tau-leaping approximates the reaction firing counts by independent Poisson random variables,}
\be
\met{\X(t+\tau)=\X(t)+\sum_{\mu=1}^{M}\mathbf{s}_{\mu}\mathcal{P}_{\mu}\!\left(a_{\mu}(\X(t))\tau\right),}
\ee
\met{where $\mathcal{P}_{\mu}(\lambda)$ is the Poisson random variable with mean $\lambda$; see \cite{gillespie2001approximate,cao2006efficient}. We use tau-leaping in two different ways. The first is adaptive tau-leaping, where the internal step size is chosen automatically to satisfy the leaping condition. This is used in Table \ref{tab:time} to compare the number of simulation steps. The second is fixed-step tau-leaping with $\tau=\Delta$, which is used as an equal-step baseline against sFML to compare approximation accuracy in Table \ref{tab:global}. In this paper, the tau-leaping method is implemented with the Python package GillesPy2 \texttt{TauLeapingSolver} \cite{Matthew2023}.}

\met{Table \ref{tab:time} compares the simulation steps of SSA, adaptive tau-leaping, and sFML. The SSA and adaptive tau-leaping steps are averaged over simulated trajectories, while the sFML steps are fixed by $T/\Delta$. Across the examples, sFML takes substantially fewer steps than SSA and adaptive tau-leaping.}

\begin{table}[htbp]
\centering
\caption{Comparison of simulation steps. SSA and tau-leaping steps are estimated as averages over 10,000 simulated trajectories and rounded to two significant figures. sFML steps equal $T/\Delta$ exactly. The ratios for SSA and adaptive tau-leaping are defined as their respective steps divided by the sFML steps.}
\label{tab:time}
\begin{tabular}{lcccccc}
\toprule
\multirowcell{2}{Example} 
& \multirowcell{2}{$T$} 
& \multirowcell{2}{sFML Steps} 
& \multicolumn{2}{c}{SSA} 
& \multicolumn{2}{c}{Adaptive Tau-leaping} \\
\cmidrule(lr){4-5} \cmidrule(lr){6-7}
& & & Steps & Ratio & Steps & Ratio \\
\midrule
Transfer & $10$ & $1.0 \times 10^{2}$ & $1.9 \times 10^{2}$ & $1.93\times$ & $1.8 \times 10^{2}$ & $1.75\times$ \\
LV (slow) & $20$ & $2.0 \times 10^{2}$ & $7.2 \times 10^{3}$ & $36\times$ & $2.5 \times 10^{3}$ & $12.3\times$ \\
LV (fast) & $20$ & $2.0 \times 10^{3}$ & $6.0 \times 10^{5}$ & $300\times$ & $3.2 \times 10^{4}$ & $15.8\times$ \\
Brusselator & $15$ & $1.5 \times 10^{3}$ & $1.6 \times 10^{6}$ & $1070\times$ & $2.1 \times 10^{5}$ & $143\times$ \\
Autocatalysis & $10$ & $1.0 \times 10^{3}$ & $3.2 \times 10^{5}$ & $317\times$ & $7.6 \times 10^{3}$ & $7.65\times$ \\
Oregonator & $6$ & $6.0 \times 10^{2}$ & $7.1 \times 10^{5}$ & $1187\times$ & $9.0 \times 10^{4}$ & $149\times$ \\
\bottomrule
\end{tabular}%
\end{table}

\met{We next evaluate the local accuracy of the learned one-step transition kernel. For a fixed conditioning state $\x$, let $\{\X_{1}^{(i)}\}_{i=1}^{m}$ be samples generated by SSA from $P_\Delta(\x,\cdot)$, and let $\{\widetilde\X_{1}^{(j)}\}_{j=1}^{n}$ be samples generated by $\G_\Theta(\x,\xxi)$. We compare these two empirical distributions using the empirical estimator of the (square of) Maximum Mean Discrepancy (MMD), defined by}
\be
\met{
\begin{split}
\widehat{\mathrm{MMD}}_h^2
=&\frac{1}{m^2}\sum_{i,i'=1}^{m} k_h\!\left(\X_{1}^{(i)},\X_{1}^{(i')}\right)
+\frac{1}{n^2}\sum_{j,j'=1}^{n} k_h\!\left(\widetilde\X_{1}^{(j)},\widetilde\X_{1}^{(j')}\right)\\
&-\frac{2}{mn}\sum_{i=1}^{m}\sum_{j=1}^{n} k_h\!\left(\X_{1}^{(i)},\widetilde\X_{1}^{(j)}\right),
\end{split}
}
\ee
\met{with the Gaussian kernel}
\be
\met{k_h(\x,\y)=\exp\left(-\frac{\norm{\x-\y}_2^2}{2h^2}\right).}
\ee
\met{The number reported in Table \ref{tab:onestep} is $\widehat{\mathrm{MMD}}_h$, where $h$ is set to the median pairwise distance of the samples \cite{gretton2012kernel}. For each example, we take one SSA test trajectory, randomly select 30 conditioning states along that trajectory, and compute the MMD at each selected state with $10,000$ conditional samples for both SSA and sFML. Table \ref{tab:onestep} reports the minimum, median, and maximum MMD values over these 30 conditioning states. These values reflect the one-step distributional accuracy along the test manifolds sampled by the SSA trajectories.}

\begin{table}[htbp]
\centering
\caption{One-step transition-kernel comparison. For each example, MMD is computed at 30 conditioning states sampled from one SSA test trajectory. The table reports the minimum, median, and maximum MMD over the selected states.}
\label{tab:onestep}
\begin{tabular}{lccc}
\toprule
Example & Min MMD & Median MMD & Max MMD \\
\midrule
Transfer & $1.89 \times 10^{-2}$ & $5.59 \times 10^{-2}$ & $1.05 \times 10^{-1}$ \\
LV (slow) & $1.72 \times 10^{-2}$ & $4.46 \times 10^{-2}$ & $9.60 \times 10^{-2}$ \\
LV (fast) & $1.55 \times 10^{-2}$ & $5.07 \times 10^{-2}$ & $8.72 \times 10^{-2}$ \\
Brusselator & $8.09 \times 10^{-2}$ & $2.29 \times 10^{-1}$ & $2.83 \times 10^{-1}$ \\
Autocatalysis & $1.81 \times 10^{-2}$ & $4.18 \times 10^{-2}$ & $9.35 \times 10^{-2}$ \\
Oregonator & $5.28 \times 10^{-2}$ & $1.53 \times 10^{-1}$ & $3.09 \times 10^{-1}$ \\
\bottomrule
\end{tabular}
\end{table}

\met{Finally, we compare long-time ensemble statistics over the full recursive rollout. Let $t_k=k\Delta$, $k=0,\ldots,K$, and define}
\be
\met{\m_k=\mathbb{E}(\X_k),\qquad \widetilde{\m}_k=\mathbb{E}(\widetilde\X_k),\qquad \s_k=\mathrm{STD}(\X_k),\qquad \widetilde{\s}_k=\mathrm{STD}(\widetilde\X_k).}
\ee
\met{The trajectory-level relative errors are}
\be
\met{E_\mu=\left(\frac{\sum_{k=0}^{K}\norm{\m_k-\widetilde{\m}_k}_2^2}{\sum_{k=0}^{K}\norm{\m_k}_2^2}\right)^{1/2},\qquad E_\sigma=\left(\frac{\sum_{k=0}^{K}\norm{\s_k-\widetilde{\s}_k}_2^2}{\sum_{k=0}^{K}\norm{\s_k}_2^2}\right)^{1/2}.}
\ee
\met{These errors measure the discrepancy of the mean and standard-deviation curves over the whole discrete trajectory. Table \ref{tab:global} reports these errors for sFML and fixed-step tau-leaping with the same coarse step size $\tau=\Delta$. This provides a quantitative counterpart to the mean and standard-deviation plots shown below for each example.}

\begin{table}[htbp]
\centering
\caption{Global trajectory comparison with fixed-step tau-leaping. Tau-leaping uses the same coarse step size $\tau=\Delta$ as sFML. Errors are estimated from 10,000 trajectories.}
\label{tab:global}
\begin{tabular}{lcccrr}
\toprule
Example & $T$ & $\Delta$ & Method & $E_\mu$ & $E_\sigma$ \\
\midrule
\multirow{2}{*}{Transfer}
& \multirow{2}{*}{$10$} 
& \multirow{2}{*}{$0.1$} 
& sFML & $1.78 \times 10^{-3}$ & $7.92 \times 10^{-2}$ \\
& & & Tau-leaping & $7.64 \times 10^{-3}$ & $4.81 \times 10^{-2}$ \\
\midrule
\multirow{2}{*}{LV (slow)}
& \multirow{2}{*}{$20$} 
& \multirow{2}{*}{$0.1$} 
& sFML & $2.16 \times 10^{-2}$ & $3.07 \times 10^{-2}$ \\
& & & Tau-leaping & $1.70 \times 10^{-1}$ & $3.04 \times 10^{0}$ \\
\midrule
\multirow{2}{*}{LV (fast)}
& \multirow{2}{*}{$20$} 
& \multirow{2}{*}{$0.01$} 
& sFML & $8.31 \times 10^{-3}$ & $5.43 \times 10^{-2}$ \\
& & & Tau-leaping & unstable & unstable \\
\midrule
\multirow{2}{*}{Brusselator}
& \multirow{2}{*}{$15$} 
& \multirow{2}{*}{$0.01$} 
& sFML & $6.87 \times 10^{-2}$ & $9.31 \times 10^{-2}$ \\
& & & Tau-leaping & $2.44 \times 10^{-1}$ & $2.93 \times 10^{-1}$ \\
\midrule
\multirow{2}{*}{Autocatalysis}
& \multirow{2}{*}{$10$} 
& \multirow{2}{*}{$0.01$} 
& sFML & $2.18 \times 10^{-2}$ & $6.91 \times 10^{-2}$ \\
& & & Tau-leaping & $2.68 \times 10^{-1}$ & $3.67 \times 10^{0}$ \\
\midrule
\multirow{2}{*}{Oregonator}
& \multirow{2}{*}{$6$} 
& \multirow{2}{*}{$0.01$} 
& sFML & $4.41 \times 10^{-1}$ & $2.78 \times 10^{-1}$ \\
& & & Tau-leaping & $3.72 \times 10^{-1}$ & $3.92 \times 10^{-1}$ \\
\bottomrule
\end{tabular}%
\end{table}

\subsection{Transfer Process}\label{ex:transfer}
We first consider the following transfer process
\begin{subequations}
    \begin{align}
        & X_1 \xrightarrow{c_1} X_2 \\
        & X_2 \xrightarrow{c_2} X_3,
    \end{align}
\end{subequations}
where there are 3 species and 2 reactions, with reaction rates $c_1=c_2=1.0$. This system has a conserved total molecule count $X_1+X_2+X_3$, so the dimension-reduction procedure of Section \ref{compdetail} is applied. A total of 40,000 training samples are generated with initial conditions drawn from $\mathcal{U}([0,100]\times [0,60] \times [50,180])$. The time step is 
set to $\Delta=0.1$, and the model is trained for $200,000$ epochs.

\begin{figure}[htbp]
  \centering
  \label{fig:Ex1}
  \includegraphics[width=.48\textwidth]{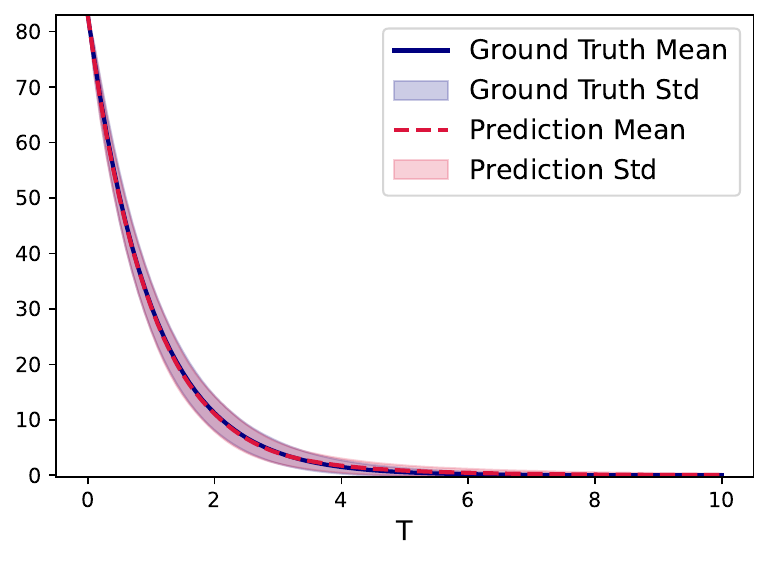}
  \includegraphics[width=.48\textwidth]{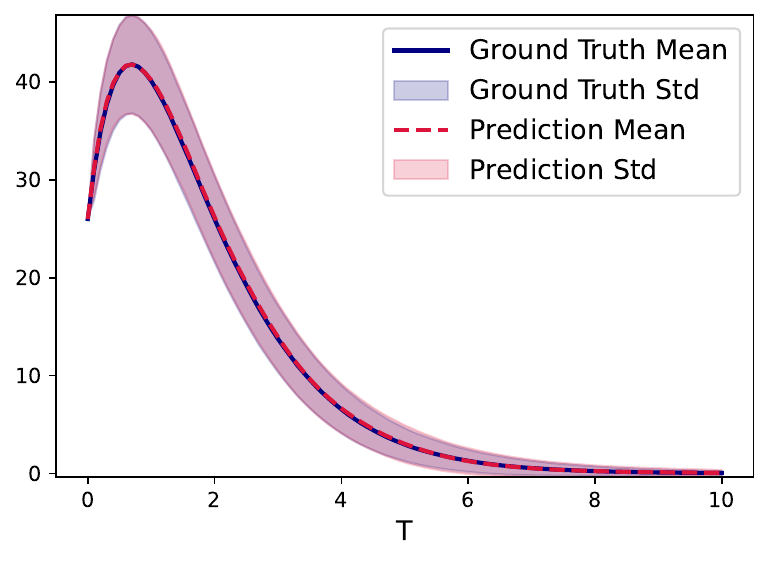}
  \includegraphics[width=.48\textwidth]{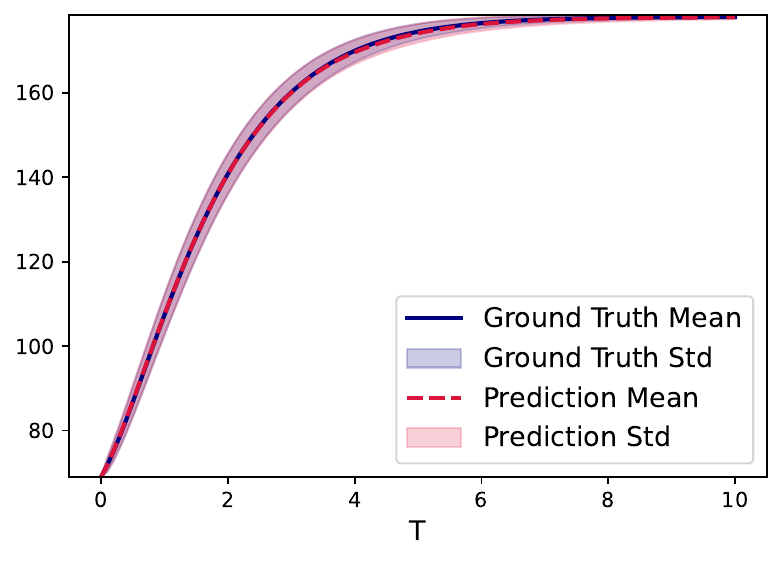}
  \includegraphics[width=.48\textwidth]{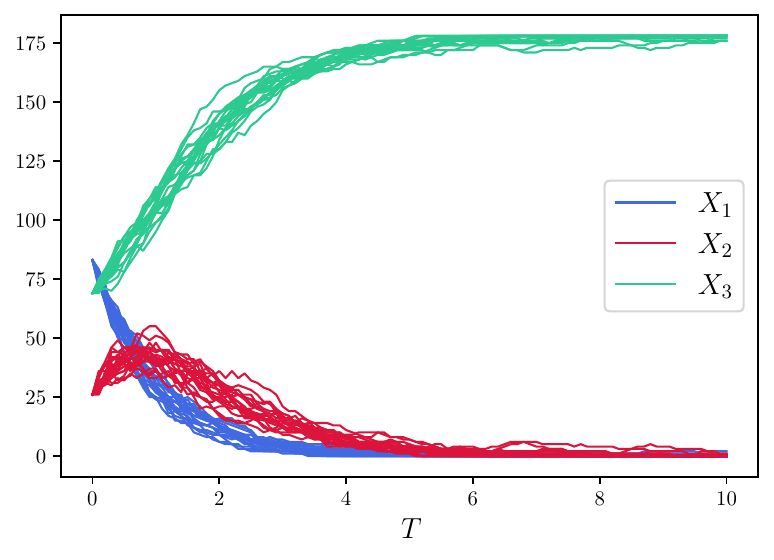}
  \caption{Mean and standard deviation (STD) of Example \ref{ex:transfer} with initial condition $(83,26,69)$: upper left: $X_1$, upper right: $X_2$, lower left: $X_3$. Lower right: sample trajectories of Example \ref{ex:transfer} with initial condition $(83,26,69)$.}
\end{figure}

The trained model is tested with the initial condition $(83,26,69)$. The lower right panel of Figure 
\ref{fig:Ex1} shows several independent predicted trajectory samples up to 
$T=10.0$, along with the mean and standard deviation averaged over $10,000$ 
trajectories. The predictions are in close agreement with the reference 
statistics over a long time period, which demonstrates the accuracy of our proposed method.

\subsection{Lotka-Volterra Model}\label{ex:LV}
We next consider the Lotka-Volterra model
\begin{subequations}
    \begin{align}
        X_1     & \xrightarrow{c_1} 2 X_1 \\
        X_1+X_2 & \xrightarrow{c_2} 2 X_2 \\
        X_2     & \xrightarrow{c_3} \emptyset,
    \end{align}
\end{subequations}
where there are 2 species and 3 reactions with rates $c_1$, $c_2$, and $c_3$. This is the classical biological model for the interaction between a predator population ($X_2$) and a prey population ($X_1$). We examine two parameter regimes that pose qualitatively different challenges to the learned model.

\begin{figure}[htbp]
  \centering
  \label{fig:Ex2_1_s}
  \includegraphics[width=.7\textwidth]{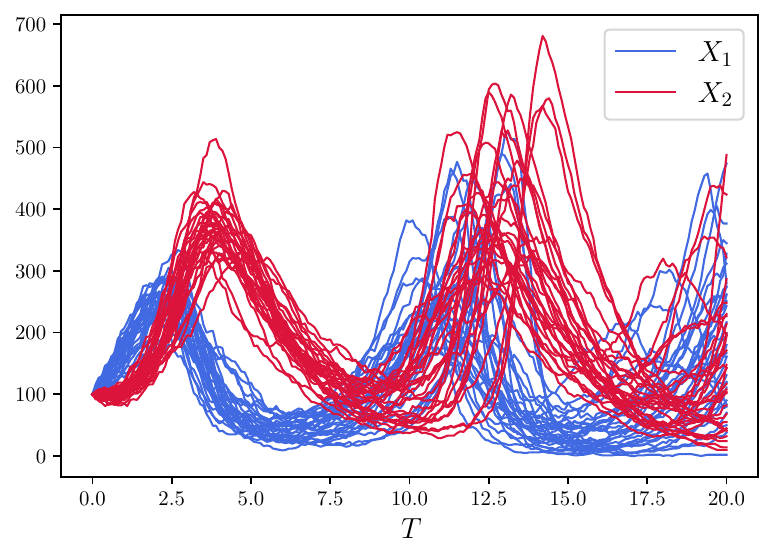}
  \caption{Test sample trajectories of Example \ref{ex:LV} with initial condition $(100,100)$ and parameter $c_1=1.0$, $c_2=0.005$ and $c_3=0.6$ up to $T=20$.}
\end{figure}

\begin{figure}[htbp]
  \centering
  \label{fig:Ex2_1}
  \includegraphics[width=.48\textwidth]{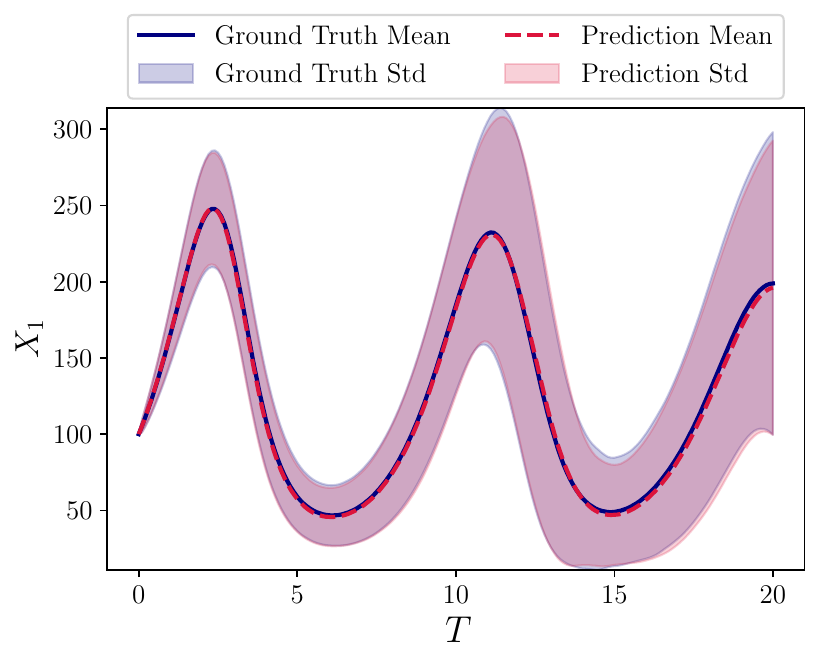}
  \includegraphics[width=.48\textwidth]{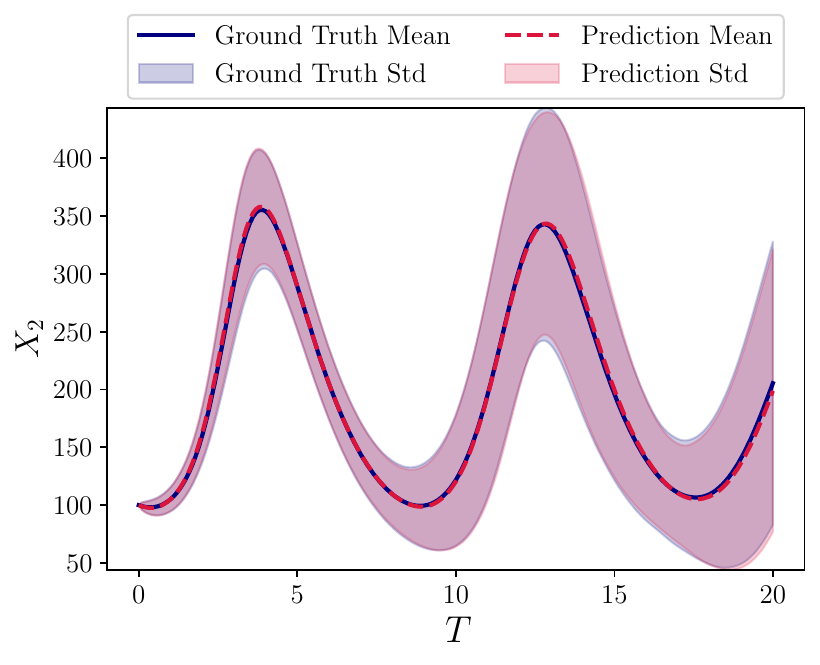}
  \caption{Mean and standard deviation (STD) of Example \ref{ex:LV} with initial condition $(100,100)$ and parameters $c_1=1.0$, $c_2=0.005$ and $c_3=0.6$ up to $T=20$: left: number of $X_1$, right: number of $X_2$.}
\end{figure}

\subsubsection*{Slow Reaction}
With $c_1=1.0$, $c_2=0.005$, and $c_3=0.6$, the system evolves on a long time scale with mild, irregular oscillations. A total of 120,000 training 
samples are generated with initial conditions drawn from 
$\mathcal{U}([0,650]\times [0,700])$. The time step is $\Delta=0.1$, and 
the model is trained for $200,000$ epochs.

The trained model is tested with the initial condition $\X_0=(100,100)$. Figure 
\ref{fig:Ex2_1_s} shows independent sample trajectories produced by the learned model, 
and Figure \ref{fig:Ex2_1} presents the mean and standard deviation 
comparison using those samples. The results are in close agreement with the reference.

\subsubsection*{Fast Reaction}
We consider the parameter set from \cite{gillespie1977exact}, with $c_1=10$, 
$c_2=0.01$, and $c_3=10$. Under such a setting, the system exhibits sharp transients and rapid oscillations. A total of 120,000 training samples are generated 
with initial conditions drawn from $\mathcal{U}([0,5000]^2)$. The time step 
is $\Delta=0.01$, chosen to resolve the sharp transients in the dynamics, 
and the model is trained for $400,000$ epochs.

\begin{figure}[htbp]
  \centering
  \label{fig:Ex2_2}
  \includegraphics[width=.48\textwidth]{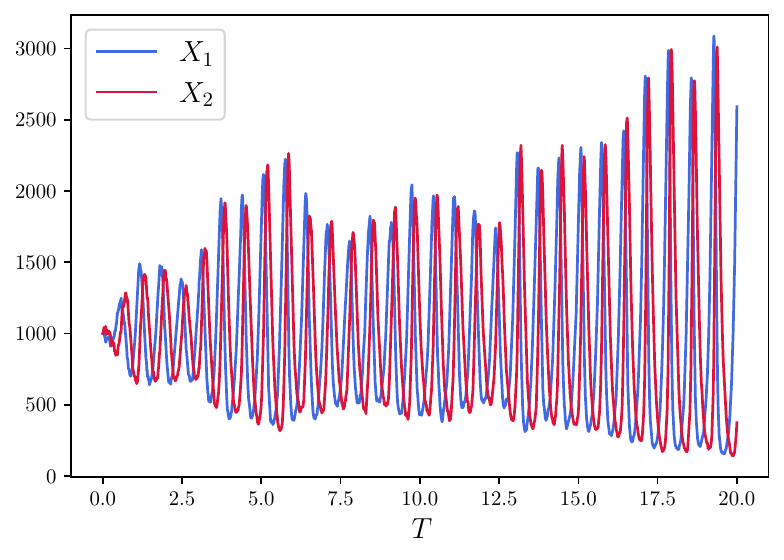}
  \includegraphics[width=.48\textwidth]{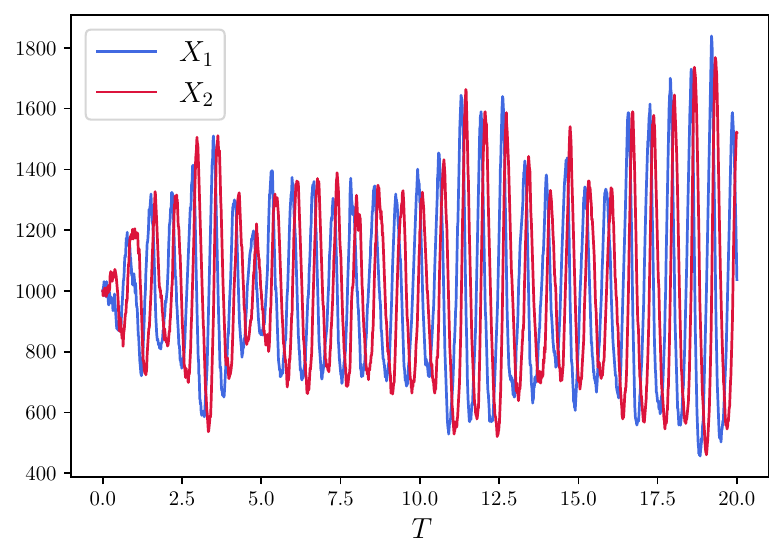}
  \includegraphics[width=.48\textwidth]{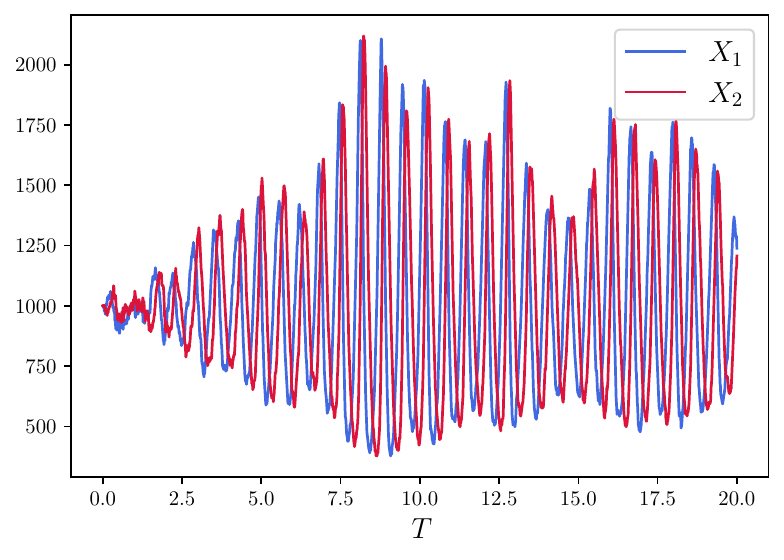}
  \includegraphics[width=.45\textwidth]{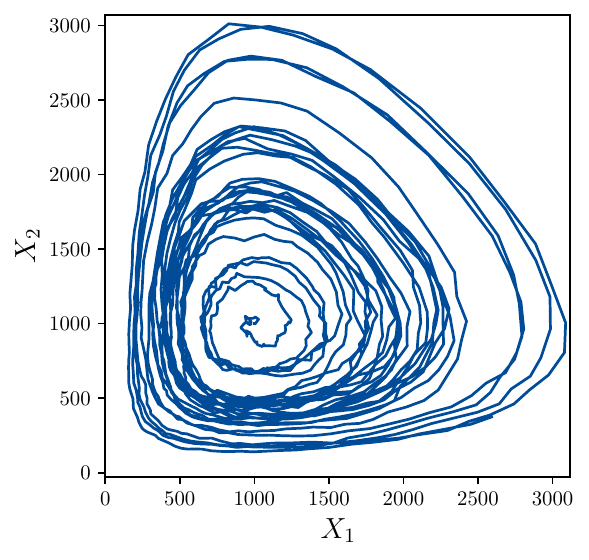}
  \caption{Sample trajectories produced by trained model of Example \ref{ex:LV} with initial condition $(1000,1000)$ and parameter $c_1=10$, $c_2=0.01$ and $c_3=10$ up to $T=20$. Upper left, upper right, lower left: 3 independent simulations, lower right: phase plot of upper left.}
\end{figure}

\begin{figure}[htbp]
  \centering
  \label{fig:Ex2_2_g}
  \includegraphics[width=.48\textwidth]{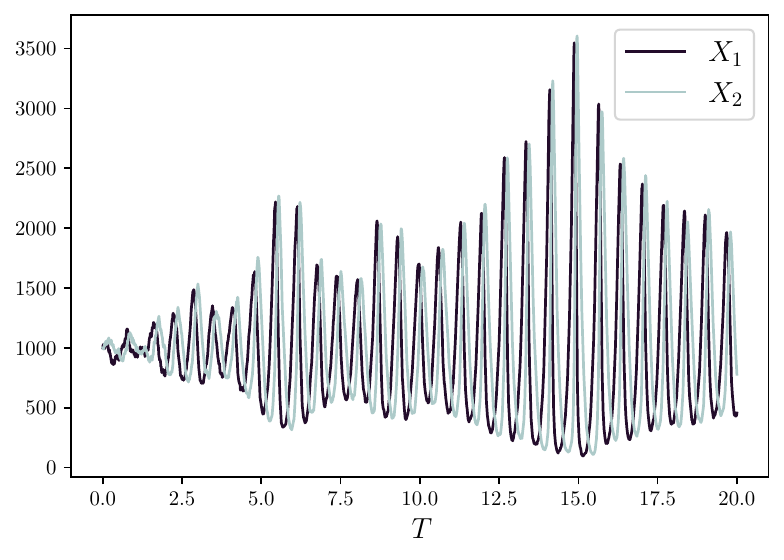}
  \includegraphics[width=.48\textwidth]{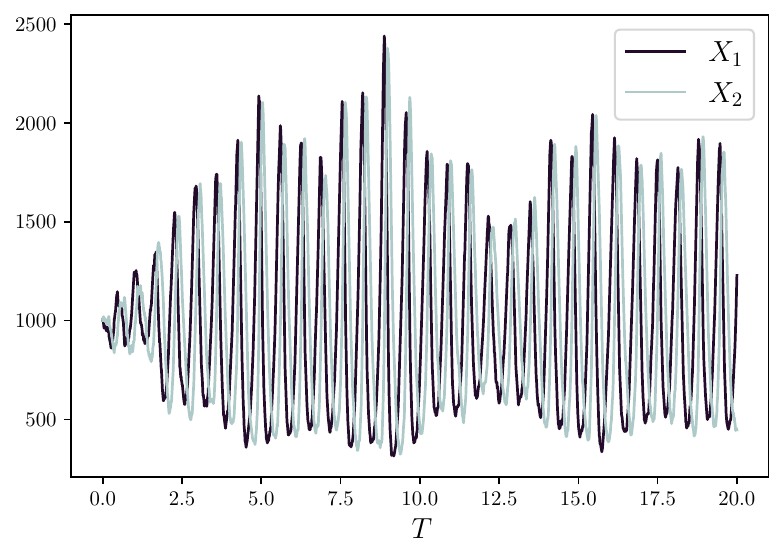}
  \includegraphics[width=.48\textwidth]{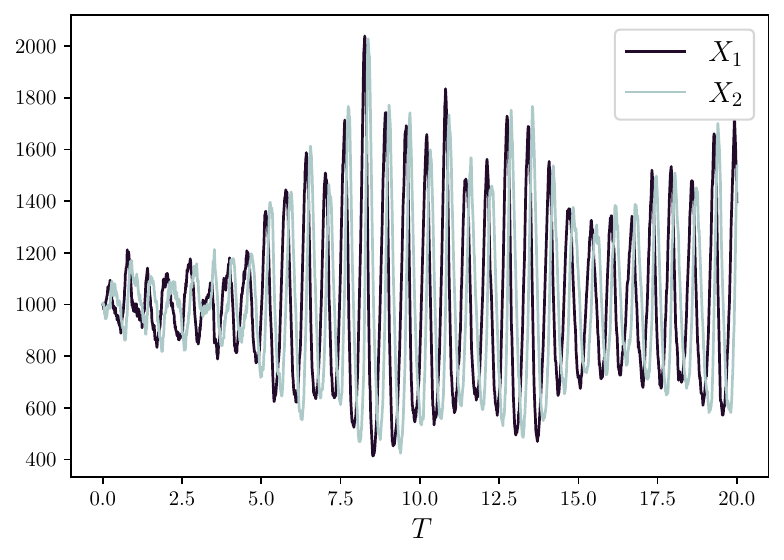}
  \includegraphics[width=.45\textwidth]{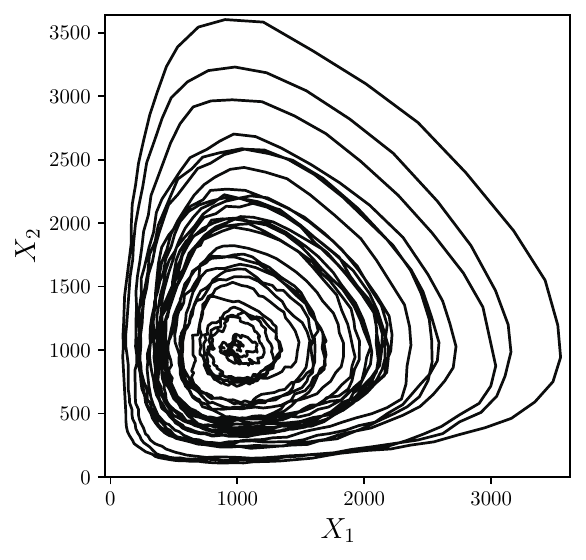}
  \caption{Sample trajectories produced by SSA of Example \ref{ex:LV} with initial condition $(1000,1000)$ and parameter $c_1=10$, $c_2=0.01$ and $c_3=10$ up to $T=20$. Upper left, upper right, lower left: 3 independent simulations, lower right: phase plot of upper left.}
\end{figure}

\begin{figure}[htbp]
  \centering
  \label{fig:Ex2_2_cond}
  \includegraphics[width=.48\textwidth]{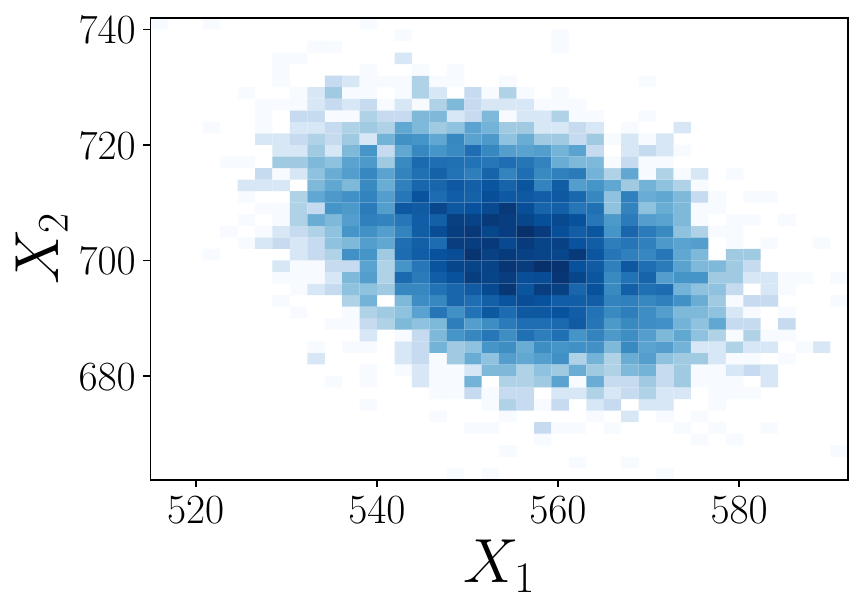}
  \includegraphics[width=.48\textwidth]{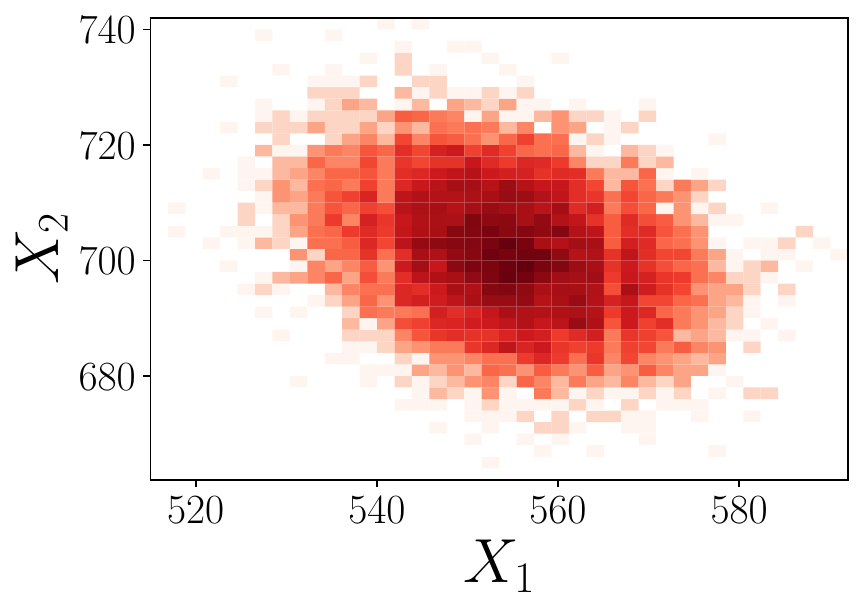}
  \includegraphics[width=.48\textwidth]{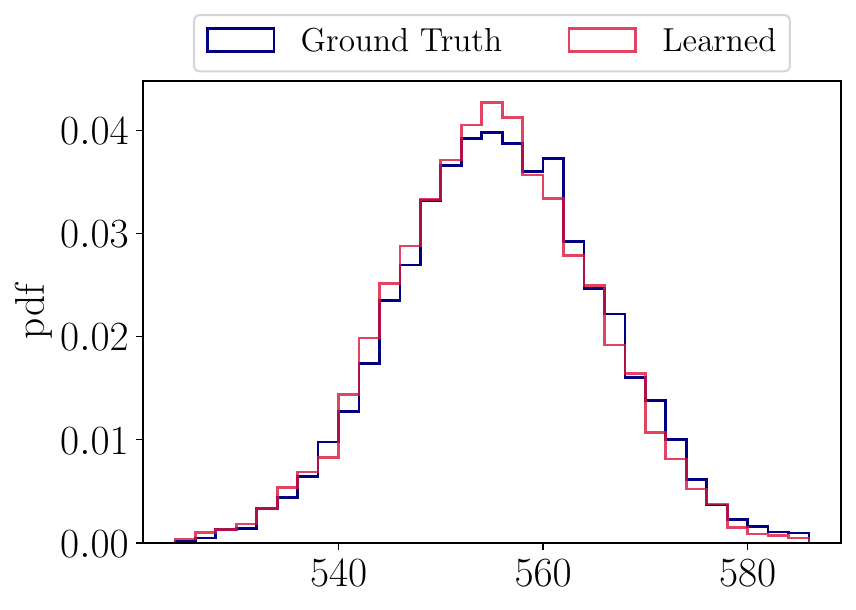}
  \includegraphics[width=.48\textwidth]{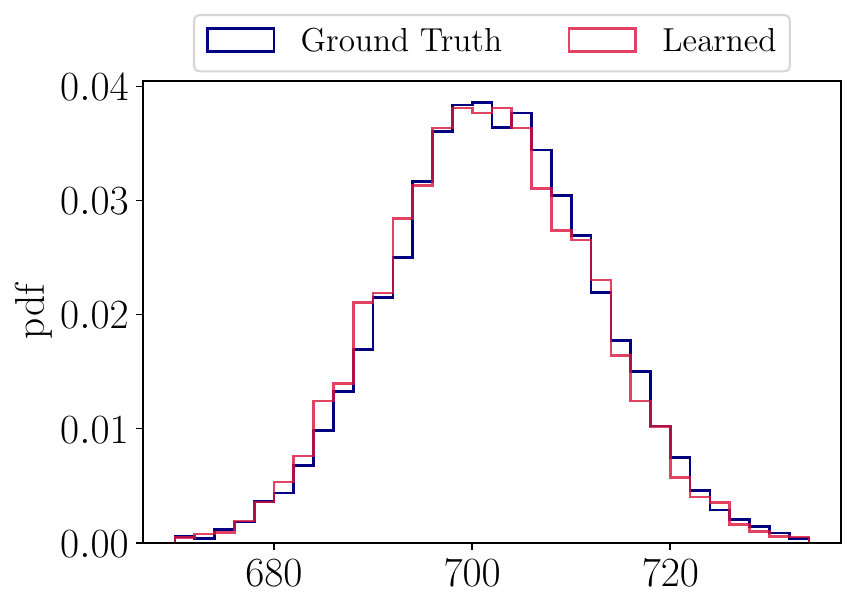}
  \caption{Comparison of one-step conditional distribution of Example \ref{ex:LV} with condition $\x=(540,734)$ and parameter $c_1=10$, $c_2=0.01$ and $c_3=10$. Upper left: ground truth of joint distribution of $(X_1,X_2)$, upper right: learned joint distribution of $(X_1,X_2)$, lower left: comparison of the marginal distribution for $X_1$, lower right: comparison of the marginal distribution for $X_2$.}
\end{figure}

Upon successful training, the model is tested with an initial condition $\X_0=(1000,1000)$. Figure \ref{fig:Ex2_2} presents 3 independent trajectory samples in 
state space and the phase space plot of the upper left case produced by the learned model. Figure \ref{fig:Ex2_2_g} shows the corresponding SSA reference trajectories. The 
two sets exhibit consistent qualitative behavior. To further assess accuracy, 
we compare the one-step conditional distributions at condition $\x=(540,734)$: 
10,000 samples are drawn from the trained model $\G_\Theta(\x,\cdot)$ and compared against the SSA referenced transition kernel $P_\Delta(\x,\cdot)$ in Figure \ref{fig:Ex2_2_cond}. The marginal distributions of $X_1$ and $X_2$ are also shown and confirm the accuracy of the learned model.

\subsection{The Brusselator}\label{ex:bru}
We consider the Brusselator
\begin{subequations}
    \begin{align}
        \overline{Y}_1          & \xrightarrow{c_1} X_1 \\
        \overline{Y}_2 + X_1    & \xrightarrow{c_2} X_2 + Z_1 \\
        2X_1 + X_2              & \xrightarrow{c_3} 3 X_1 \\
        X_1                     & \xrightarrow{c_4} Z_2,
    \end{align}
\end{subequations}
where there are 2 species and 4 reactions, with $\overline{Y}_1$ and 
$\overline{Y}_2$ treated as parameters and $c_1$, $c_2$, $c_3$, and $c_4$ the reaction rates. A total of 120,000 training samples 
are generated with initial conditions drawn from $\mathcal{U}([0,5000]^2)$, 
parameter values $c_1\overline{Y}_1=5000$, $c_2\overline{Y}_2=50$, 
$c_3=5\times10^{-5}$, $c_4=5$. The time step is set to $\Delta=0.01$, and the model is trained for $400,000$ epochs.

\begin{figure}[htbp]
  \centering
  \label{fig:Ex3}
  \includegraphics[width=.48\textwidth]{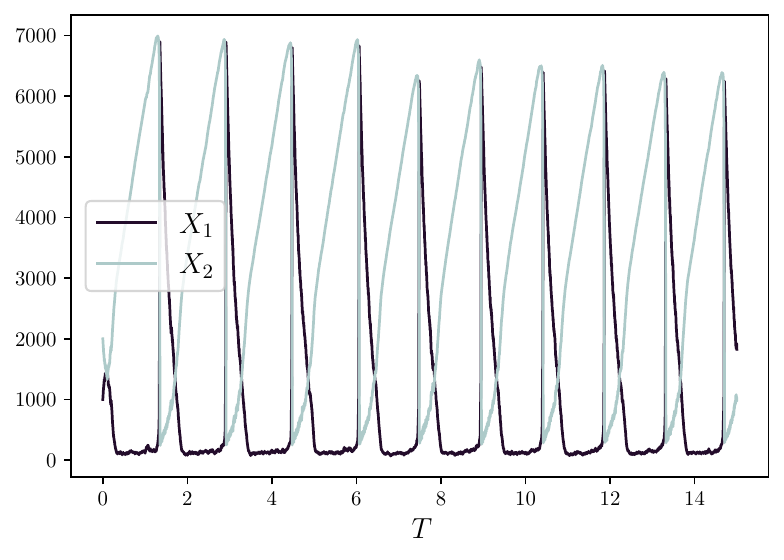}
  \includegraphics[width=.48\textwidth]{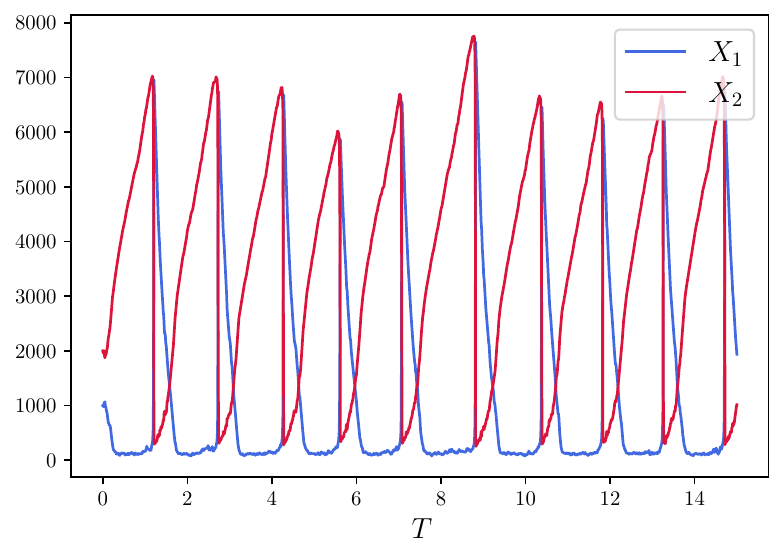}
  \includegraphics[width=.48\textwidth]{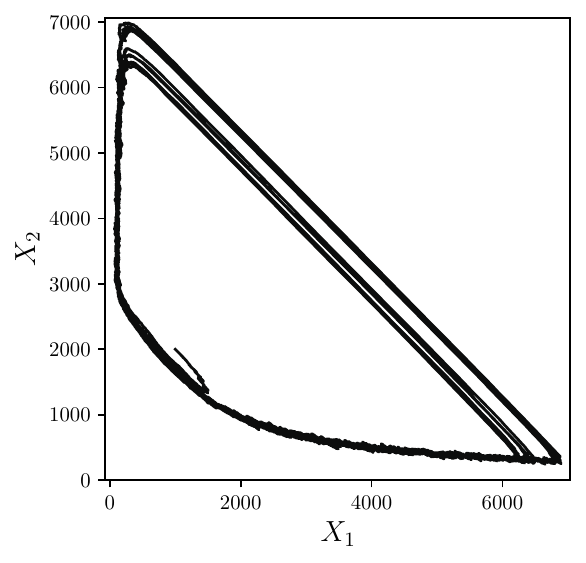}
  \includegraphics[width=.48\textwidth]{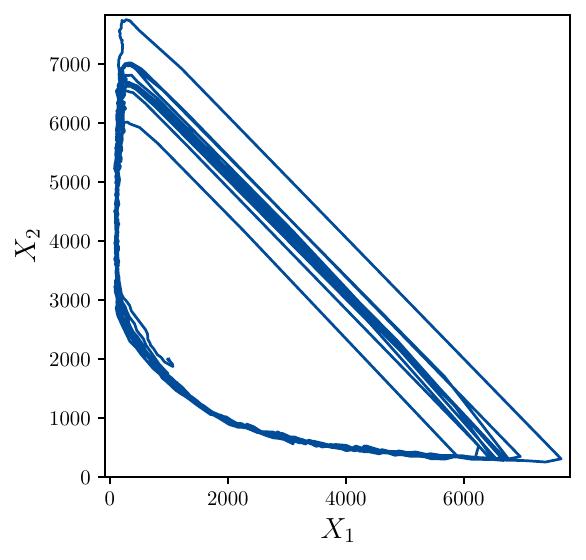}
  \caption{Comparison of sample trajectories of Example \ref{ex:bru} with initial condition $(1000,2000)$ and parameter $c_1 \overline{Y}_1=5000$, $c_2 \overline{Y}_2=50$, $c_3=5\times 10^{-5}$ and $c_4=5$ up to $T=15$. Upper left: a sample trajectory in state space produced by SSA, upper right: a sample trajectory in state space produced by learned model, lower left: phase plot of upper left, lower right: phase plot of upper right.}
\end{figure}

\begin{figure}[htbp]
  \centering
  \label{fig:Ex3_cond}
  \includegraphics[width=.48\textwidth]{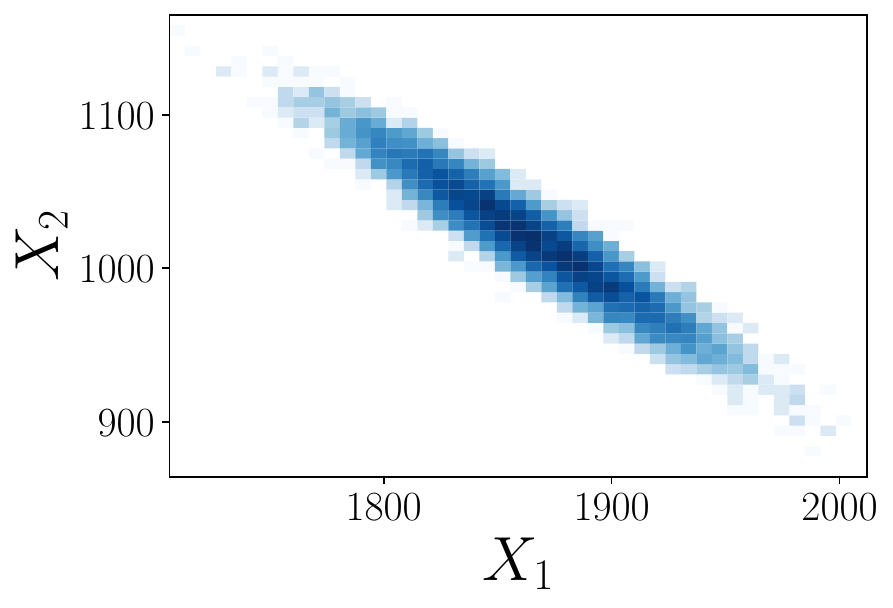}
  \includegraphics[width=.48\textwidth]{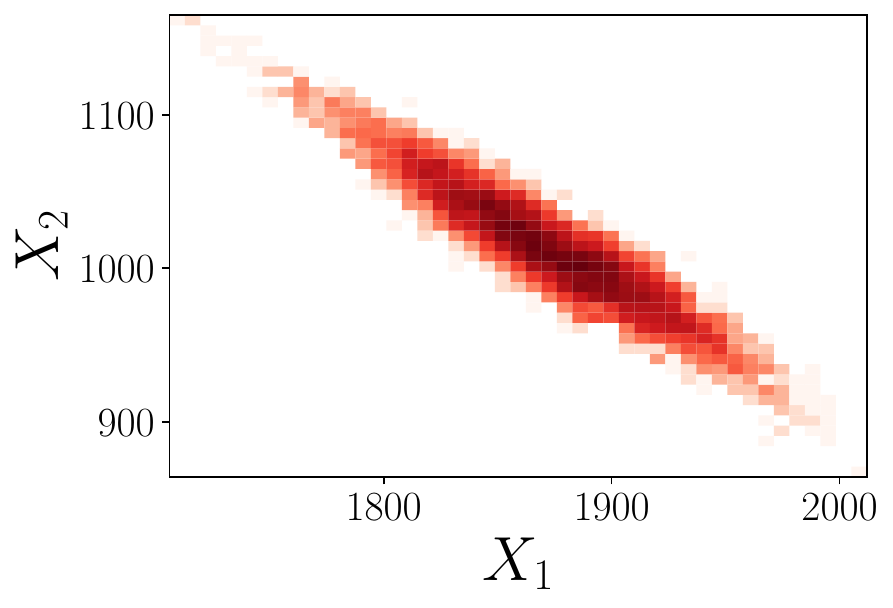}
  \includegraphics[width=.48\textwidth]{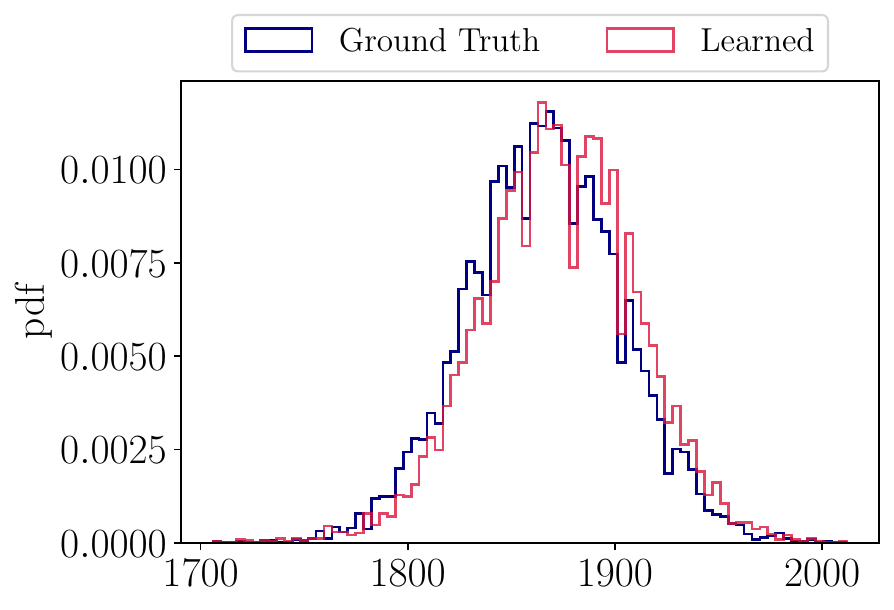}
  \includegraphics[width=.48\textwidth]{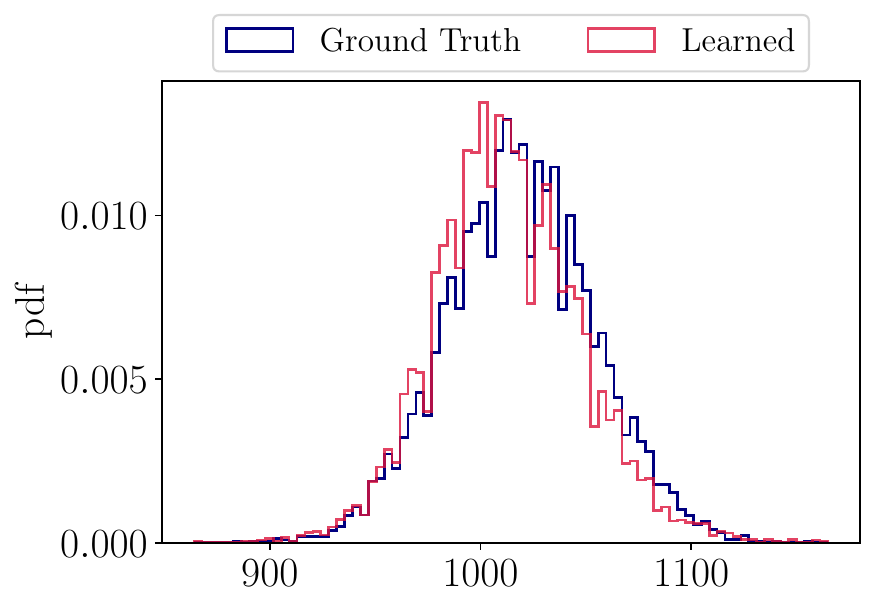}
  \caption{Comparison of one-step conditional distribution of Example \ref{ex:bru} with initial condition $(1959,971)$ and parameter $c_1 \overline{Y}_1=5000$, $c_2 \overline{Y}_2=50$, $c_3=5\times 10^{-5}$ and $c_4=5$. Upper left: ground truth of joint distribution of $(X_1,X_2)$, upper right: learned joint distribution of $(X_1,X_2)$, lower left: comparison of the marginal distribution for $X_1$, lower right: comparison of the marginal distribution for $X_2$.}
\end{figure}

\begin{figure}[htbp]
  \centering
  \label{fig:Ex3_fft}
  \includegraphics[width=.48\textwidth]{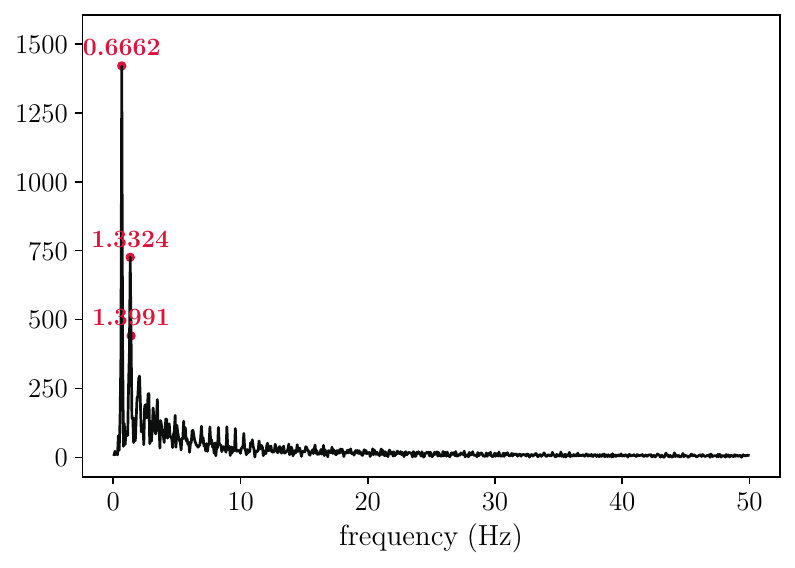}
  \includegraphics[width=.48\textwidth]{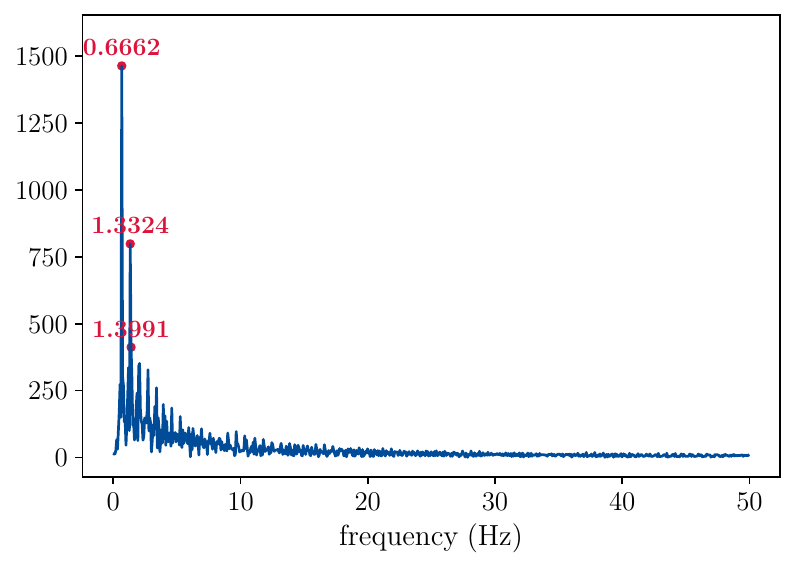}
  \caption{Comparison of Fourier modes for sFML generated trajectory and ground truth of Example \ref{ex:bru} with initial condition $(1000,2000)$ and parameter $c_1 \overline{Y}_1=5000$, $c_2 \overline{Y}_2=50$, $c_3=5\times 10^{-5}$ and $c_4=5$ up to $T=15$. Left: Fourier modes of SSA trajectory (upper left of Figure \ref{fig:Ex3}), right: Fourier modes of model trajectory (upper right of Figure \ref{fig:Ex3}).}
\end{figure}

Upon successful training, the model is tested with an initial condition $\X_0=(1000,2000)$. Figure \ref{fig:Ex3} presents a comparison of trajectories (in both state and phase spaces) generated by SSA and by the trained model. A consistent, nearly periodic behavior is observed across both sets of results. To further validate the dynamical properties, we apply the Fast Fourier Transform (FFT) to compute the Fourier modes of the two trajectories, presented in Figure \ref{fig:Ex3_fft}. The first, second, and third dominant modes of the two trajectories agree closely, demonstrating that the learned model has captured the correct oscillation frequency. Furthermore, a comparison of histograms of the conditional distributions is presented in Figure \ref{fig:Ex3_cond}. With condition $\x=(1959,971)$, $10,000$ samples are generated from the trained model $\G_\Theta(\x,\cdot)$ and displayed in the upper right of Figure \ref{fig:Ex3_cond}, while the corresponding ground-truth histogram obtained from SSA transition kernel $P_\Delta(\x,\cdot)$ is shown in the upper left.  The marginal distributions of $X_1$ and $X_2$ are also compared, confirming the accuracy of the learned model.

\subsection{The Autocatalysis}\label{ex:auto}
We consider the autocatalysis process
\begin{subequations}
    \begin{align}
        X_1+X_2    & \xrightarrow{c_1} 2 X_2 \\
        X_2+X_3    & \xrightarrow{c_2} 2 X_3 \\
        X_3+X_1    & \xrightarrow{c_3} 2 X_1,
    \end{align}
\end{subequations}
where there are 3 species and 3 reactions with rates $c_1=c_2=c_3=0.002$. A total of 120,000 training samples are generated with initial conditions drawn from $\mathcal{U}([1000,4000]^3)$. The time step is taken as $\Delta=0.01$. The training epochs are taken to be $400,000$.

\begin{figure}[htbp]
  \centering
  \label{fig:Ex4}
  \includegraphics[width=.48\textwidth]{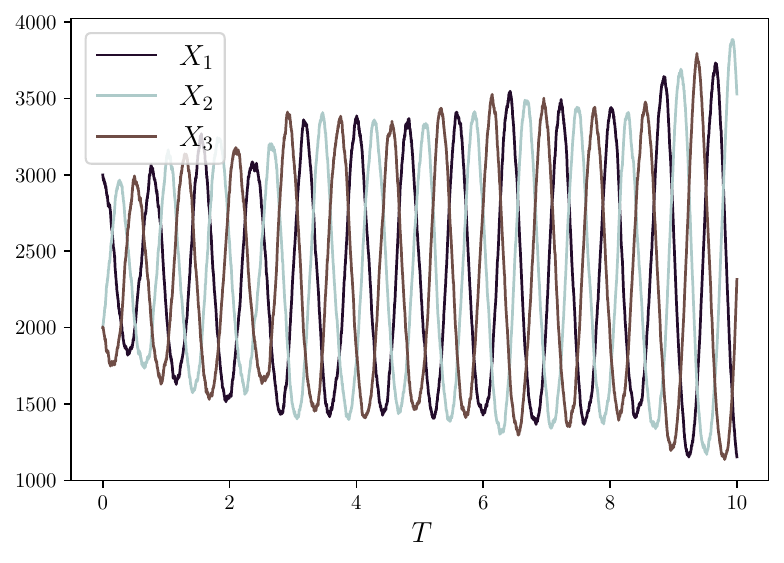}
  \includegraphics[width=.48\textwidth]{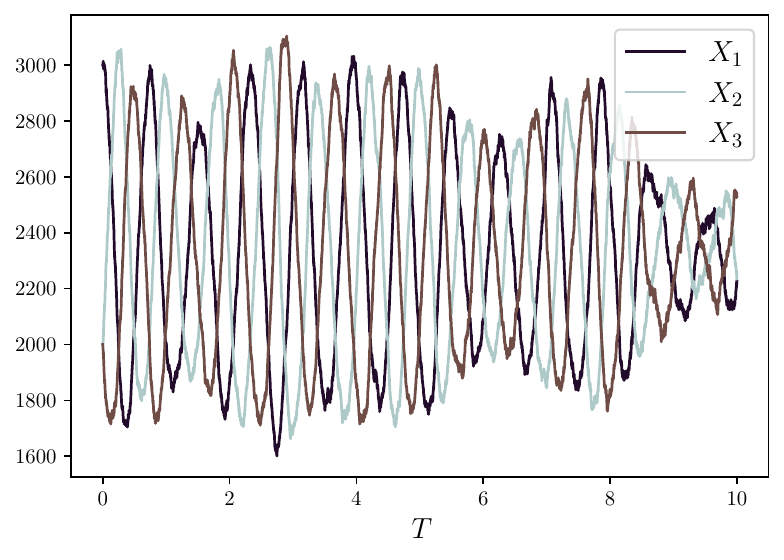}
  \includegraphics[width=.48\textwidth]{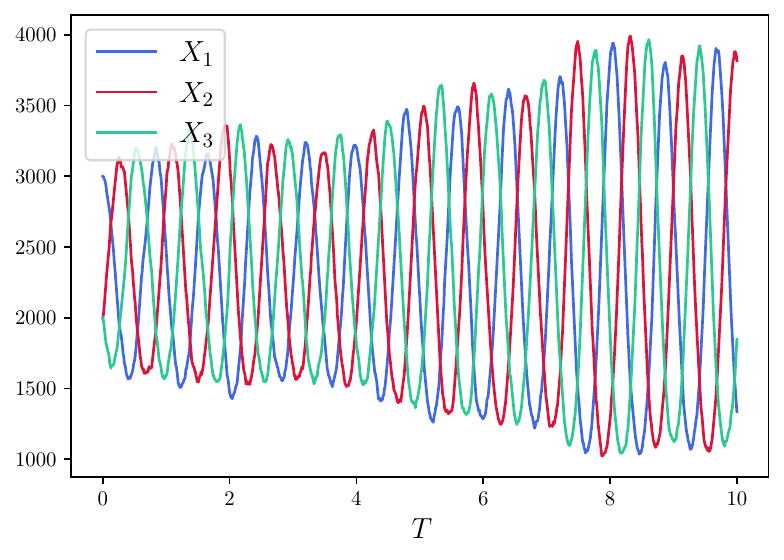}
  \includegraphics[width=.48\textwidth]{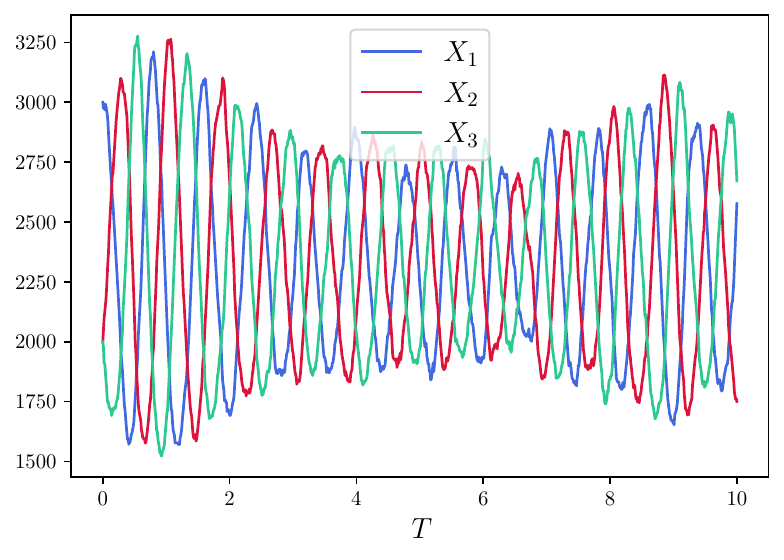}
  \caption{Comparison of sample trajectories of Example \ref{ex:auto} with initial condition $(3000,2000,2000)$ and parameter $c_1=0.002$, $c_2=0.002$ and $c_3=0.002$ up to $T=10$. Upper left, upper right: 2 independent sample trajectories produced by SSA, lower left, lower right: 2 independent sample trajectories produced by the trained model.}
\end{figure}

\begin{figure}[htbp]
  \centering
  \label{fig:Ex4_fft}
  \includegraphics[width=.48\textwidth]{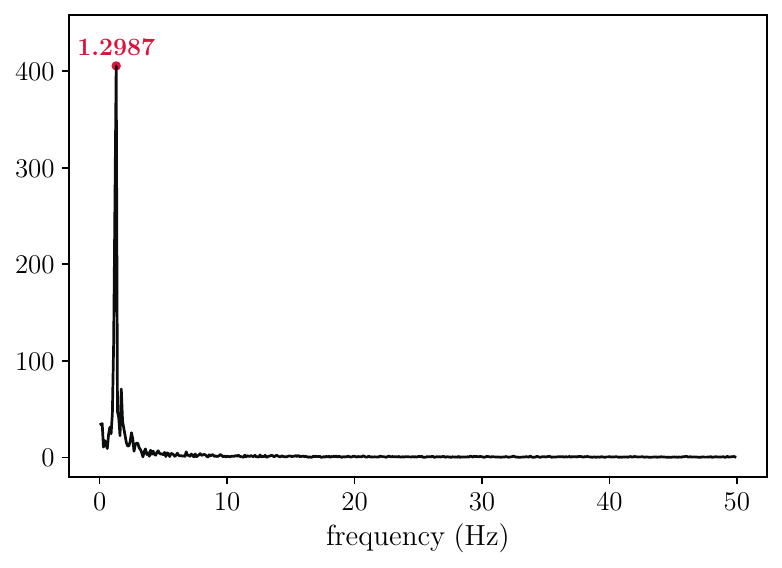}
  \includegraphics[width=.48\textwidth]{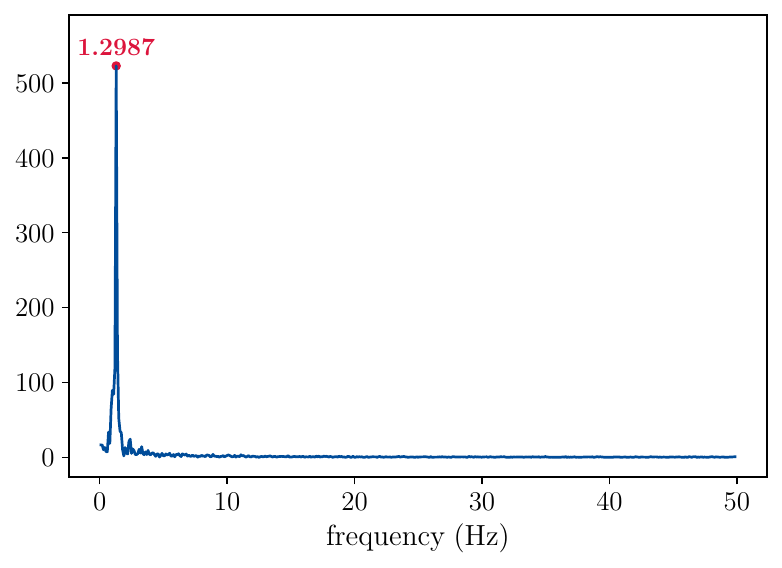}
  \caption{Comparison of Fourier modes for sFML generated trajectory and ground truth of Example \ref{ex:auto} with initial condition $(3000,2000,2000)$ and parameter $c_1=0.002$, $c_2=0.002$ and $c_3=0.002$ up to $T=10$. Left: Fourier modes of SSA trajectory (upper right of Figure \ref{fig:Ex4}), right: Fourier modes of model trajectory (lower right of Figure \ref{fig:Ex4}).}
\end{figure}

This system conserves the total molecule count, so the dimension-reduction procedure of Section \ref{compdetail} applies. The trained model is tested with an initial condition $\X_0=(3000,2000,2000)$. 
Figure \ref{fig:Ex4} compares trajectories in the state space produced 
by SSA and the learned model. The system exhibits periodic oscillations, and 
the dominant Fourier modes of the two trajectories, presented in Figure 
\ref{fig:Ex4_fft}, are in close agreement. This confirms that the cyclic oscillation frequency is correctly learned.

\subsection{The Oregonator}\label{ex:org}
We consider the Oregonator model
\begin{subequations}
    \begin{align}
        \overline{Y}_1 +X_2         & \xrightarrow{c_1} X_1 \\
        X_1 + X_2                   & \xrightarrow{c_2} Z_1 \\
        \overline{Y}_2 +X_1         & \xrightarrow{c_3} 2 X_1 +X_3 \\
        2 X_1                       & \xrightarrow{c_4} Z_2 \\
        \overline{Y}_3 + X_3        & \xrightarrow{c_5} X_2,
    \end{align}
\end{subequations}
where there are 3 species and 5 reactions, $ \overline{Y}_1$, $\overline{Y}_2$ and $\overline{Y}_3$ are treated as parameters, $c_1$, $c_2$, $c_3$, $c_4$ and $c_5$ are reaction rates, with $c_1 \overline{Y}_1=2$, $c_2 = 0.1$, $c_3 \overline{Y}_2=104$, $c_4=0.016$ and $c_5 \overline{Y}_3=26$. A total of 120,000 training samples are generated 
with initial conditions drawn from $\mathcal{U}([0,10000]^3)$. The time step is taken as $\Delta=0.01$. The training epochs are taken to be $400,000$.

\begin{figure}[htbp]
  \centering
  \label{fig:Ex5}
  \includegraphics[width=.48\textwidth]{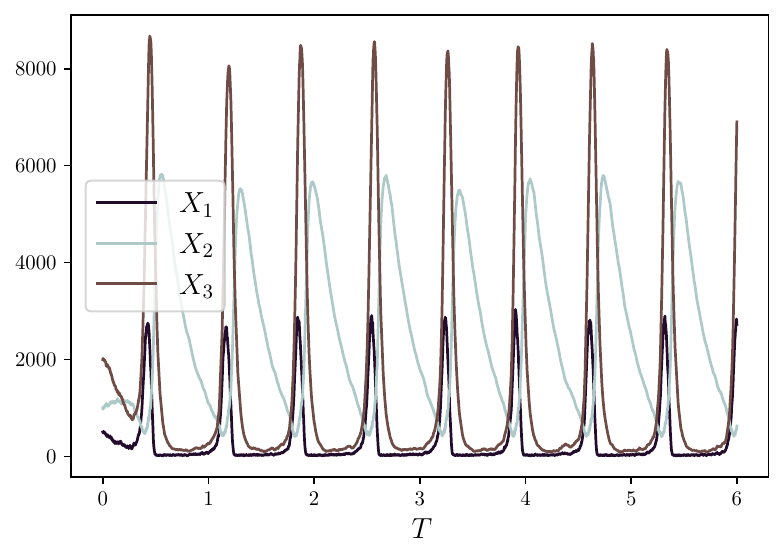}
  \includegraphics[width=.48\textwidth]{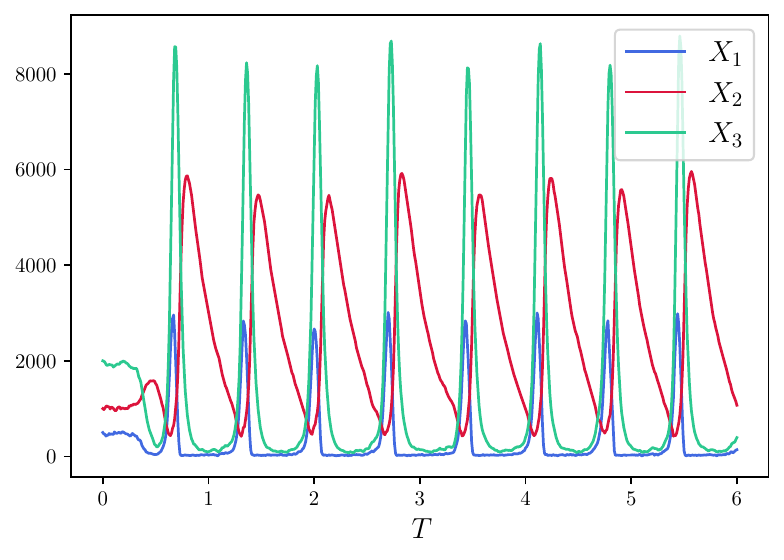}
  \caption{Comparison of sample trajectories of Example \ref{ex:org} in state space with initial condition $(500,1000,2000)$ and parameter $c_1 \overline{Y}_1=2$, $c_2 = 0.1$, $c_3 \overline{Y}_2=104$, $c_4=0.016$ and $c_5 \overline{Y}_3=26$ up to $T=6$. Left: a sample trajectory produced by SSA, right: a sample trajectory produced by the trained model.}
\end{figure}

\begin{figure}[htbp]
  \centering
  \label{fig:Ex5_phase}
  \includegraphics[width=.31\textwidth]{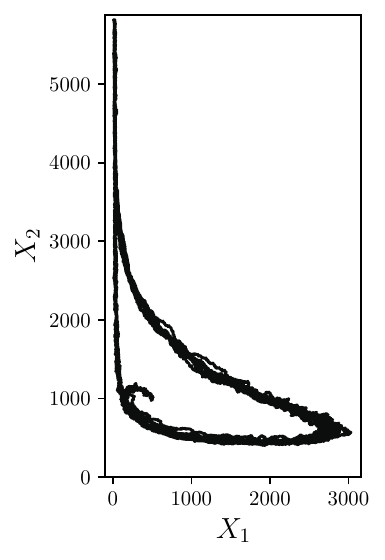}
  \includegraphics[width=.268\textwidth]{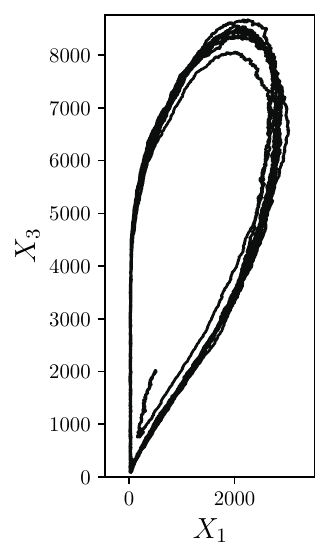}
  \includegraphics[width=.38\textwidth]{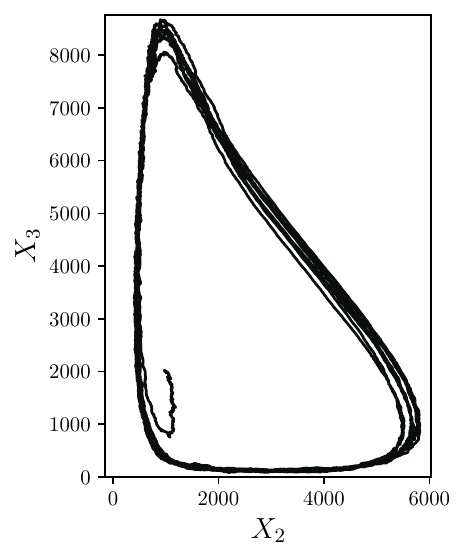}
  \includegraphics[width=.31\textwidth]{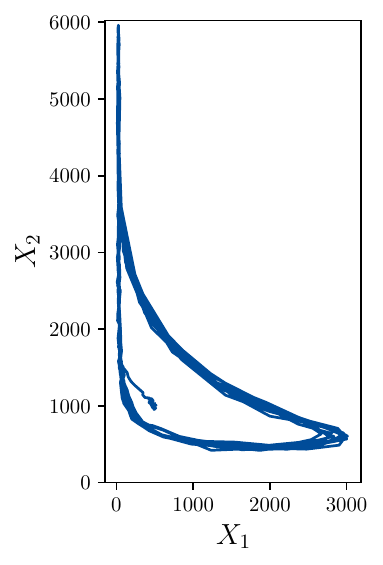}
  \includegraphics[width=.268\textwidth]{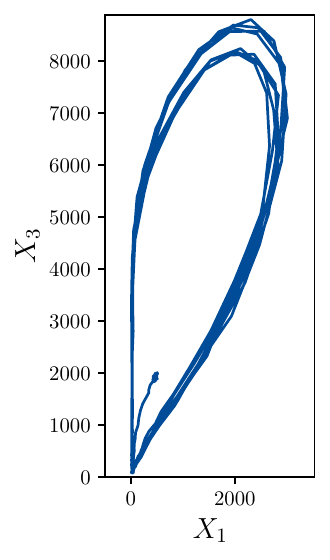}
  \includegraphics[width=.38\textwidth]{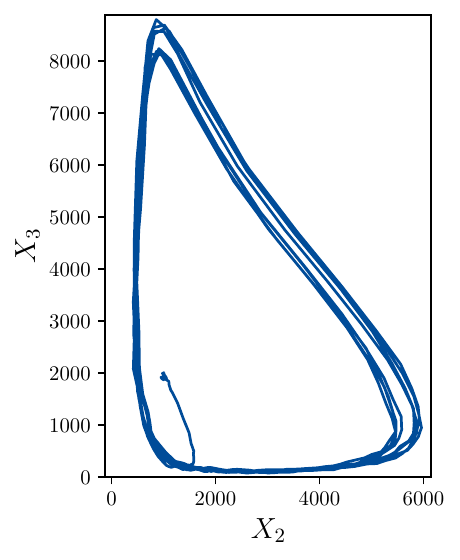}
  \caption{Comparison of sample trajectories of Example \ref{ex:org} in phase space with initial condition $(500,1000,2000)$ and parameter $c_1 \overline{Y}_1=2$, $c_2 = 0.1$, $c_3 \overline{Y}_2=104$, $c_4=0.016$ and $c_5 \overline{Y}_3=26$ up to $T=6$. Upper: sample trajectory produced by SSA, lower: sample trajectory produced by trained model.}
\end{figure}

\begin{figure}[htbp]
  \centering
  \label{fig:Ex5_fft}
  \includegraphics[width=.48\textwidth]{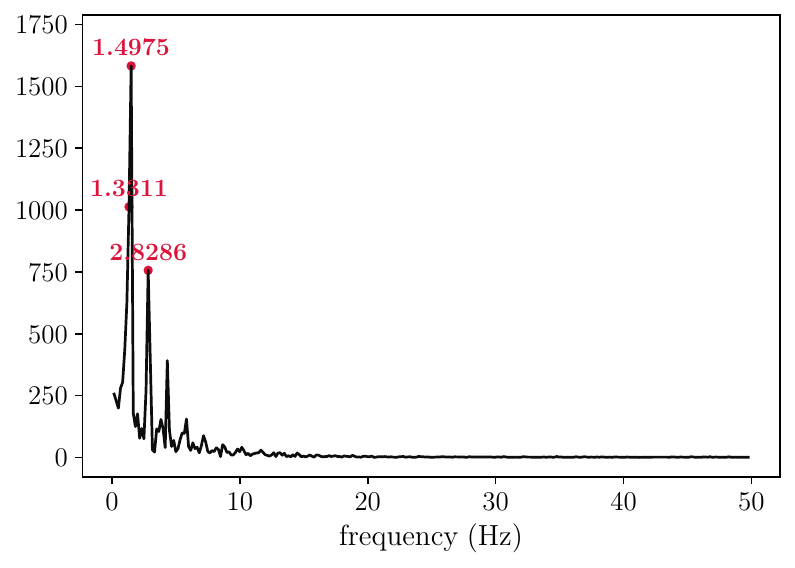}
  \includegraphics[width=.48\textwidth]{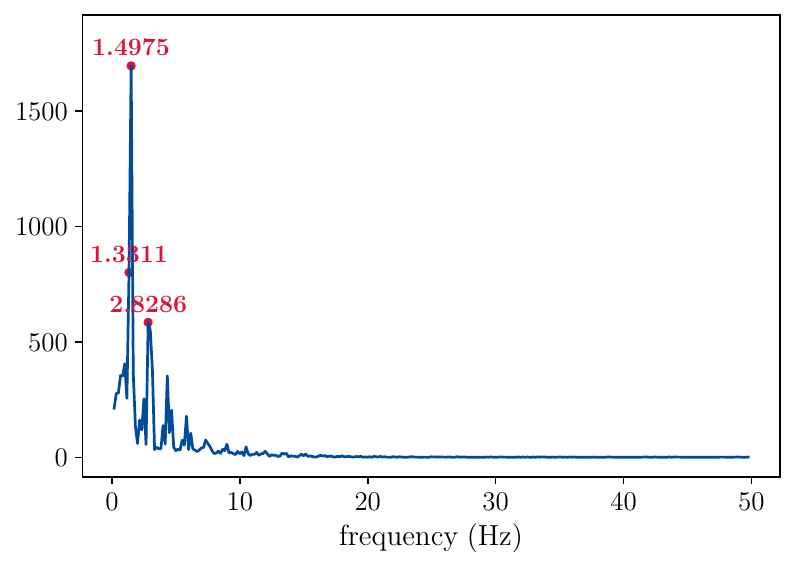}
  \caption{Comparison of Fourier modes for two sample trajectories of Example \ref{ex:org} with initial condition $(500,1000,2000)$ and parameter $c_1 \overline{Y}_1=2$, $c_2 = 0.1$, $c_3 \overline{Y}_2=104$, $c_4=0.016$ and $c_5 \overline{Y}_3=26$ up to $T=6$. Left: Fourier modes of SSA trajectory (left of Figure \ref{fig:Ex5}), right: Fourier modes of model trajectory (right of Figure \ref{fig:Ex5}).}
\end{figure}

The trained model is tested with initial condition $\X_0=(500,1000,2000)$. In Figures \ref{fig:Ex5} and \ref{fig:Ex5_phase}, we show a comparison of trajectories (in state and phase spaces) with both SSA and our trained model. The dominant Fourier modes of the two trajectories, presented in Figure \ref{fig:Ex5_fft}, agree closely, confirming the accuracy of the learned model.

\section{Conclusion} \label{sec:conclu}
In this paper, we propose a data-driven effective model for the Stochastic Simulation Algorithm of stochastic reaction networks. The method directly approximates the transition kernel of the underlying continuous-time Markov chain from short bursts of SSA simulation data using a generative model. This approach naturally accommodates highly non-Gaussian transition distributions and provides a flexible framework for constructing coarse stochastic models from exact SSA data. The resulting model serves as a stochastic propagator that generates long-time trajectories at a user-defined, constant coarse time step, independent of the microscopic firing times. As a consequence, the number of simulation steps is reduced significantly while the statistical properties of the true dynamics are accurately reproduced. In this paper, we employ a conditional normalizing flow as the realization of such a stochastic propagator. By using several numerical examples, we demonstrate the effectiveness and efficiency of the proposed method.

\section*{Data Availability Statement}
The data that support the findings of this study are available from the corresponding author upon reasonable request.

\bibliographystyle{siamplain}
\bibliography{references}
\end{document}